\documentclass[a4paper,fleqn]{cas-sc}

\usepackage{graphicx}
\usepackage{amsmath,amssymb,amsfonts}
\usepackage{amsthm}
\usepackage{algorithm}
\usepackage{algorithmic}
\usepackage{mathtools}
\usepackage{braket}
\usepackage{xcolor}
\usepackage{textcomp}
\usepackage{booktabs}
\usepackage[numbers,sort&compress]{natbib}
\usepackage{placeins}

\theoremstyle{definition}

\newcommand{\p}{\partial}
\renewcommand{\d}{\mathrm{d}}
\renewcommand{\O}[1]{\mathcal{O}\left(#1\right)}
\newcommand{\dS}{\mathrm{dS}}
\newcommand{\x}{\mathbf{x}}

\newcommand{\io}{\int_\Omega}
\newcommand{\bs}[1]{\boldsymbol{#1}}
\newcommand{\wt}[1]{\widetilde{#1}}
\newcommand{\T}{\mathrm{T}}
\newcommand{\GRobin}{\Gamma^{\mathrm{Robin}}}
\newcommand{\GDirichlet}{\Gamma^{\mathrm{Dirichlet}}}
\newcommand{\uhom}{\hat{u}^{\mathrm{hom}}}
\newcommand{\ulift}{\hat{u}^{\mathrm{lift}}}
\newcommand{\Ahom}{A^{\mathrm{hom}}}
\newcommand{\bhom}{b^{\mathrm{hom}}}

\begin{document}

\shorttitle{Tensor-Network Chebyshev--Galerkin Method}
\shortauthors{Rodríguez-Aldavero et al.}

\title[mode=title]{A Tensor-Network Chebyshev--Galerkin Method for Boundary-Value Problems in Computational Mechanics}

\author[1,2]{Juan José Rodríguez-Aldavero}[orcid=0009-0001-8942-4285]
\cormark[1]
\ead{jj.aldavero@csic.es}

\author[3]{Ángel de la Rosa Velasco}

\author[1,2]{Juan José García-Ripoll}

\affiliation[1]{organization={Quantum Advanced Research Center (QuARC), CSIC},
  addressline={Calle Serrano 113b},
  city={Madrid},
  postcode={28006},
  country={Spain}}

\affiliation[2]{organization={Institute of Fundamental Physics (IFF), CSIC},
  addressline={Calle Serrano 113b},
  city={Madrid},
  postcode={28006},
  country={Spain}}

\affiliation[3]{organization={Grupo de Durabilidad e Integridad Mecánica de Materiales Estructurales (DIMME), Universidad Rey Juan Carlos},
  addressline={Calle Tulipán s/n},
  city={Móstoles},
  postcode={28933},
  country={Spain}}

\cortext[1]{Corresponding author}

\begin{abstract}
Quantum-inspired numerical methods for partial differential equations represent discretized fields and operators as quantics tensor trains (QTT), so that their cost depends on how compressible the problem is rather than on the number of unknowns. Global spectral formulations in this setting have so far relied on Fourier bases, restricted to periodic problems, or on collocation, which does not preserve the variational structure exploited by tensor-network solvers. We introduce a QTT Chebyshev--Galerkin method that assembles and solves global spectral systems directly in compressed form, without forming the underlying arrays. Starting from the weak formulation, the method handles coupled vector-valued elliptic equations with variable coefficients, smoothly mapped geometries, and mixed Dirichlet--Robin boundary conditions. Boundary conditions are imposed through algebraic operations on compressed operators rather than entrywise modifications, preserving the symmetric positive-definite structure exploited by variational tensor-network solvers. Parameters and random inputs enter as additional coordinates, so the same construction returns the solution over an entire parameter space in a single solve. In linear elasticity, the method retains spectral convergence and reaches errors below $10^{-8}$ on a three-dimensional benchmark, beyond those attained by QTT finite differences and Q2/Q4 FEniCSx finite elements over the tested range. A six-parameter cantilever design study covering $1.7\times10^{7}$ designs is solved in a single 92~s run, and a stochastic problem with ten random coordinates is represented with $2.4\times10^{6}$ times less memory than its full coefficient array.
\end{abstract}


\maketitle

\section{Introduction}
\label{sec:introduction}

Linear elasticity describes the displacements, strains, and stresses induced in a solid by body forces, surface tractions, and prescribed displacements~\cite{slaughter2002}, and is a prototypical elliptic boundary-value problem (BVP) of computational mechanics~\cite{evans2022}. Its numerical solution becomes demanding in two ways. First, accurate solutions may require fine resolution, particularly in three dimensions, where material properties or loads vary over length scales much smaller than the body, as in functionally graded materials~\cite{birman2007,jha2013}. Second, many applications require the solution not for a single configuration but over a range of inputs: mechanical design explores geometric, material, and boundary parameters~\cite{ion2022ttiga}, while uncertainty quantification propagates random material properties, such as a spatially varying elastic modulus~\cite{ghanem2003,matthies2005}. Treating each parameter as an additional coordinate, the solution lives on a joint spatial--parametric domain, whose tensor-product discretization with $N$ points along each of $d$ coordinates contains $N^d$ unknowns, a number that grows rapidly with resolution and exponentially with the number of coordinates. The conventional alternative, typically Monte Carlo sampling, solves one configuration at a time, but its error decreases only algebraically with the number of samples, so accurate results require a very large number of solves~\cite{caflisch1998}.

Finite-difference (FD) and finite-element (FEM) methods approximate the solution locally, leading to sparse linear systems whose error decreases algebraically with resolution~\cite{zienkiewicz1977,hughes2003}; high-order and spectral elements raise this rate within a local discretization~\cite{patera1984,karniadakis2005}. Global spectral methods instead expand the solution in smooth basis functions spanning the whole domain, such as Fourier modes or Chebyshev polynomials~\cite{trefethen2000,boyd2001}. For smooth solutions on simple or smoothly mapped domains, their error decreases exponentially with the number of basis functions, so high accuracy requires far fewer unknowns per direction. Spectral Galerkin methods add a variational formulation that incorporates boundary conditions naturally and, as in FEM, yields symmetric positive-definite systems~\cite{shen1995,shen2011}. Fast convergence, however, reduces the number of unknowns per direction but not their exponential growth with the number of coordinates, and, since every basis function spans the whole domain, the resulting operators are dense: a three-dimensional elasticity problem with 100 modes per direction already involves $3\times10^6$ coefficients and a dense operator with about $10^{13}$ entries. Structured and matrix-free solvers avoid storing such operators~\cite{olver2013,kronbichler2019}, but the solution itself remains exponentially large, which rules out high resolutions and joint spatial--parametric domains.

Tensor networks provide compressed representations of large numerical arrays that exhibit low-rank structure~\cite{kolda2009,orus2014}. In this work, they are used to mitigate the rapid growth of the coefficient arrays arising from spectral discretizations. Tensor trains (TT)~\cite{oseledets2011}, also known in quantum physics as matrix-product states (MPS)~\cite{verstraete2008,schollwock2011}, represent an array as a chain of small tensors whose size, called the TT rank or bond dimension, measures the correlations between coordinates. Storage and computational cost are then dictated by the bond dimension rather than by the nominal size of the array. Quantics tensor trains (QTT) further split each coordinate into binary digits~\cite{oseledets2010,khoromskij2011b}, so that $N$ points along one direction are represented by only $\log_2 N$ small tensors. For bounded bond dimension, storage grows only logarithmically with resolution and linearly with the number of coordinates, instead of exponentially. Linear operators admit the same compressed format, known as matrix-product operators (MPOs), enabling linear algebra directly in compressed form up to a controlled truncation error~\cite{garcia2021,garcia2026}. A QTT over $n$ binary digits has the same structure as the state of an $n$-qubit quantum register~\cite{schon2005,malz2024}. QTT methods are therefore \emph{quantum-inspired}: they borrow tools from quantum physics and offer a natural bridge to quantum computing.

Several PDE solvers have been formulated in TT and QTT form. Most rely on local discretizations, such as finite differences~\cite{dolgov2012parabolic,kazeev2012} and finite elements~\cite{kazeev2018,kazeev2022}, recently applied to linear elasticity~\cite{benvenuti2025}. These treat boundary conditions through local modifications of stencils or elements, but inherit the algebraic convergence of the underlying schemes. Global spectral formulations, which could converge exponentially, remain comparatively scarce, and a central difficulty is the treatment of non-periodic boundary conditions. Fourier-based QTT methods avoid it by assuming periodicity~\cite{garcia2023,josien2026}, while pseudospectral and collocation approaches~\cite{adak2024,gidi2025,cavallini2026} impose boundary conditions in strong form, replacing or penalizing selected rows of the discrete operator. Such entrywise manipulations are natural for explicit matrices but awkward in compressed MPO formats, where individual entries are not directly accessible, and they generally destroy the symmetry of the problem. A systematic treatment of general boundary conditions in global spectral QTT discretizations has therefore been missing. The Galerkin formulation provides one: boundary conditions enter the weak form either as boundary integrals added to the operator or as constraints on the approximation space, and can thus be assembled like the rest of the problem, without penalties or entrywise modifications.

In this work, we introduce a QTT Chebyshev--Galerkin method that assembles and solves global spectral systems directly in compressed form, without ever forming the full arrays. The method has three main features. First, it is systematic: starting from the weak formulation, problem data, bulk operators, and boundary terms are all assembled in QTT form, for coupled vector-valued elliptic equations with variable coefficients, smoothly mapped geometries, and mixed Dirichlet--Robin boundary conditions. Second, it is structure-preserving: for symmetric coercive problems such as elasticity, the resulting linear systems are symmetric positive-definite, and their solution minimizes a quadratic energy. This allows the use of variational tensor-network algorithms, routinely used for ground-state search in computational quantum physics~\cite{white1992,dolgov2014}. Third, parameters and random inputs enter as additional coordinates, so that a single solve returns the solution over an entire parameter space, extending stochastic Galerkin methods~\cite{khoromskij2011,dolgov2015pce} to a fully spectral setting. To the best of our knowledge, this is the first systematic QTT formulation of a global spectral Galerkin method that incorporates mixed Dirichlet, Neumann, and Robin boundary conditions through compressed operator algebra, while retaining the symmetry and variational structure of the underlying elliptic problem. We realize the method in the Chebyshev basis, whose algebraic structure yields accurate and compact QTT representations of differentiation, multiplication by variable coefficients, and the transform from point values to spectral coefficients. An open-source Python implementation of the framework built on SeeMPS~\cite{garcia2026}, together with all experimental data, is publicly available on Zenodo~\cite{rodriguez2026code}.

Numerical experiments in linear elasticity support two main conclusions. First, compression does not compromise accuracy: the QTT formulation retains the spectral convergence of the underlying discretization and, on a three-dimensional benchmark, reaches errors below $10^{-8}$ that QTT finite differences and quadratic and quartic finite elements based on FEniCSx~\cite{alnaes2015,baratta2023} do not attain over the tested range---although QTT finite differences remain cheaper at moderate accuracy. Second, the cost is governed by the bond dimensions of the operator and the solution rather than by the nominal number of unknowns. This is clearly observed when parametric coordinates are included, where problems far beyond the reach of explicit arrays become tractable: a single solve captures a six-parameter cantilever design space with about $1.7\times10^{7}$ designs, and a stochastic problem with ten random coordinates is stored with $2.4\times10^{6}$ times less memory than its explicit coefficient array. Taken together, when the problem exhibits low-rank spectral structure, the method combines the exponential convergence of spectral discretizations with the exponential compression of tensor networks, yielding a discretization error that decays doubly exponentially with the number of QTT cores, within a systematic and structure-preserving framework.

\section{Problem Definition}
\label{sec:problem}

We consider second-order elliptic BVPs for a coupled vector field $u(\x)=\left(u_1(\x),\ldots,u_m(\x)\right)$ with $m$ components, such as the $m=d$ displacement components of a solid in elasticity. The problem is posed on the reference domain $\Omega=[-1,1]^d$, with coordinates $\x=(x_1,\ldots,x_d)$, which is the natural domain of Chebyshev polynomials; smoothly curved physical domains are mapped onto it through a change of coordinates (see Appendix~\ref{app:mappings}). Indices $\alpha,\beta=1,\ldots,d$ label spatial directions and $\mu,\nu=1,\ldots,m$ field components; summation over repeated indices is implied. Boundary conditions may differ between components, so for each component the boundary is split into Dirichlet and Robin parts,
\begin{equation}
    \p\Omega = \GDirichlet_\mu \cup \GRobin_\mu,
    \qquad
    \GDirichlet_\mu \cap \GRobin_\mu = \varnothing,
    \qquad
    \mu = 1, \ldots, m.
    \label{eq:boundary_partition}
\end{equation}
The value of $u_\mu$ is prescribed on $\GDirichlet_\mu$, whereas a combination of its flux and value is prescribed on $\GRobin_\mu$, which includes Neumann conditions as a special case. This allows, for instance, a face to be fixed in the normal direction while sliding freely in the tangential ones. For simplicity, both parts are assumed to be unions of complete faces of $\Omega$.

\subsection{Strong Formulation}
\label{sec:problem-strong}

We consider elliptic BVPs in divergence form,
\begin{equation}
  \left\{
  \begin{alignedat}{3}
      -\p_\alpha \left[\mathbb{C}_{\mu\alpha\nu\beta}(\x) \p_\beta u_\nu(\x)\right]
      + R_{\mu\nu}(\x)u_\nu(\x)
      &={}& f_\mu(\x), &\qquad& \x &\in \Omega,
      \\[0.75em]
      n_\alpha(\x) \mathbb{C}_{\mu\alpha\nu\beta}(\x) \p_\beta u_\nu(\x)
      + S_{\mu\nu}(\x)u_\nu(\x)
      &={}& g_\mu(\x), && \x &\in \GRobin_\mu,
      \\[0.75em]
      u_\mu(\x)
      &={}& \bar{u}_\mu(\x), && \x &\in \GDirichlet_\mu.
  \end{alignedat}
  \right.
  \label{eq:bvp_strong}
\end{equation}
Here, $\mathbb{C}_{\mu\alpha\nu\beta}(\x)$ is the coefficient tensor, $R_{\mu\nu}(\x)$ the reaction tensor, and $f_\mu(\x)$ the source term. On $\GRobin_\mu$, $n_\alpha(\x)$ is the outward unit normal, $n_\alpha \mathbb{C}_{\mu\alpha\nu\beta}\p_\beta u_\nu$ the conormal flux of component $\mu$, $S_{\mu\nu}(\x)$ the Robin coefficient, and $g_\mu(\x)$ the Robin data. Setting $S_{\mu\nu}=0$ recovers Neumann conditions. On $\GDirichlet_\mu$, $\bar{u}_\mu(\x)$ is the prescribed value.

Small-strain linear elasticity is the prototypical model chosen for the numerical experiments. In this context, $u_\mu$ is the displacement field of a solid, with $m=d$. It is usually written in terms of the Cauchy stress $\sigma_{ij}=\mathbb{C}_{ijkl}\,\p_l u_k$, with indices $(i,j,k,l)$, so that equilibrium reads $-\p_j\sigma_{ij}=f_i$ and the traction on the boundary is $\sigma_{ij}n_j$. In our notation, with indices $(\mu,\alpha,\nu,\beta)$, these are the first two lines of Eq.~\eqref{eq:bvp_strong}: the coefficient tensor is the elasticity tensor, the conormal flux is the traction, $f_\mu$ and $g_\mu$ are body forces and applied tractions, $R_{\mu\nu}$ and $S_{\mu\nu}$ model elastic foundations and supports, and $\bar{u}_\mu$ is a prescribed displacement. For an isotropic material with elastic modulus $E(\x)$ and Poisson ratio $\nu_{\mathrm P}(\x)$,
\begin{equation}
    \mathbb{C}_{\mu\alpha\nu\beta}(\x)
    = \lambda_1(\x)\, \delta_{\mu\alpha} \delta_{\nu\beta}
    + \lambda_2(\x) \left( \delta_{\mu\nu} \delta_{\alpha\beta} + \delta_{\mu\beta} \delta_{\alpha\nu} \right),
    \quad
    \lambda_1(\x) = \frac{E(\x)\,\nu_{\mathrm P}(\x)}{[1+\nu_{\mathrm P}(\x)][1-2\nu_{\mathrm P}(\x)]},
    \quad
    \lambda_2(\x) = \frac{E(\x)}{2[1+\nu_{\mathrm P}(\x)]},
    \label{eq:elasticity_principal_tensor}
\end{equation}
where $\lambda_1$ and $\lambda_2$ are the Lamé parameters of the material.

\subsection{Weak Formulation and Energy Principle}
\label{sec:problem-weak}

The method is built on the weak formulation of Eq.~\eqref{eq:bvp_strong}, in which the two types of boundary conditions enter differently: Dirichlet conditions restrict the spaces of admissible fields, whereas Robin conditions contribute boundary integrals. We define the trial and test spaces
\begin{equation}
    \mathcal{V}_{\bar{u}} \coloneqq \left\{ w\in \mathrm{H}^1(\Omega,\mathbb{R}^m) : w_\mu = \bar{u}_\mu \text{ on } \GDirichlet_\mu \right\},
    \qquad
    \mathcal{V}_{0} \coloneqq \left\{ v\in \mathrm{H}^1(\Omega,\mathbb{R}^m) : v_\mu = 0 \text{ on } \GDirichlet_\mu \right\},
    \label{eq:spaces}
\end{equation}
where $\mathrm{H}^1(\Omega,\mathbb{R}^m)$ is the space of square-integrable vector fields with square-integrable first derivatives. Trial fields thus satisfy the prescribed Dirichlet data, while test fields vanish on the Dirichlet boundary. Multiplying the bulk equation in Eq.~\eqref{eq:bvp_strong} by a test field $v\in\mathcal{V}_0$, integrating over $\Omega$, and integrating by parts yields the weak problem: find $u \in \mathcal{V}_{\bar{u}}$ such that
\begin{equation}
    a(u, v) = \ell(v) \qquad \forall v \in \mathcal{V}_0,
    \label{eq:bvp_weak}
\end{equation}
with
\begin{align}
    a(u,v) &\coloneqq \io \mathbb{C}_{\mu\alpha\nu\beta}(\x) \p_\beta u_\nu(\x) \p_\alpha v_\mu(\x) \, \d\x
    \,+\, \io R_{\mu\nu}(\x) u_\nu(\x) v_\mu(\x) \, \d\x
    \,+\, \int_{\p\Omega} S_{\mu\nu}(\x) u_\nu(\x) v_\mu(\x) \, \dS,
    \label{eq:bvp_weak_bilinear}
    \\
    \ell(v) &\coloneqq \io f_\mu(\x) v_\mu(\x) \,\d\x 
    \,+\, \int_{\p\Omega} g_\mu(\x) v_\mu(\x) \,\dS.
    \label{eq:bvp_weak_functional}
\end{align}
The boundary term produced by integration by parts vanishes on $\GDirichlet_\mu$, where $v_\mu=0$, and is replaced by the Robin condition on $\GRobin_\mu$; extending $S_{\mu\nu}$ and $g_\mu$ by zero outside $\GRobin_\mu$, both boundary integrals can be written over $\p\Omega$.

We restrict to problems whose bilinear form is symmetric and coercive. Symmetry follows from
\begin{equation}
    \mathbb{C}_{\mu\alpha\nu\beta}(\x) = \mathbb{C}_{\nu\beta\mu\alpha}(\x),
    \qquad
    R_{\mu\nu}(\x) = R_{\nu\mu}(\x),
    \qquad
    S_{\mu\nu}(\x) = S_{\nu\mu}(\x),
    \label{eq:bvp_symmetry}
\end{equation}
while coercivity requires $a(v,v)\geq c\,\|v\|^2_{\mathrm{H}^1}$ for some $c>0$ and all $v\in\mathcal{V}_0$, as in linear elasticity when the boundary conditions exclude rigid-body motions. Under these assumptions, the weak solution is the unique minimizer of a quadratic energy,
\begin{equation}
    u = \arg\min_{w\in\mathcal{V}_{\bar{u}}} \mathcal{E}[w],
    \qquad
    \mathcal{E}[w] = \frac12 a(w,w) - \ell(w),
    \label{eq:energy}
\end{equation}
which in linear elasticity is the total potential energy. 

The coefficients, loads, boundary data, and geometry may also depend on parameters $\bs{\xi}=(\xi_1,\ldots,\xi_{d_\xi})\in\Xi=[-1,1]^{d_\xi}$ that are not differentiated, such as design variables or random inputs. Rather than solving Eq.~\eqref{eq:bvp_weak} separately for each $\bs{\xi}$, we seek the solution $u(\x,\bs{\xi})$ over the joint spatial--parametric domain $\Omega\times\Xi$ by also integrating the weak form over the parameters,
\begin{equation}
    \int_\Xi a(u,v;\bs{\xi})\,\d\bs{\xi} = \int_\Xi \ell(v;\bs{\xi})\,\d\bs{\xi}
    \qquad \forall v \in \mathrm{L}^2(\Xi;\mathcal{V}_0),
    \label{eq:bvp_weak_parametric}
\end{equation}
with $u\in \mathrm{L}^2(\Xi;\mathcal{V}_{\bar{u}})$, where $a(u,v;\bs{\xi})$ and $\ell(v;\bs{\xi})$ are Eqs.~\eqref{eq:bvp_weak_bilinear}--\eqref{eq:bvp_weak_functional} evaluated with the data at $\bs{\xi}$. This is the intrusive formulation of stochastic Galerkin methods~\cite{ghanem2003,matthies2005}. If $a$ is symmetric and coercive uniformly in $\bs{\xi}$, the joint problem inherits the energy principle, with energy $\int_\Xi \mathcal{E}[w;\bs{\xi}]\,\d\bs{\xi}$. Since derivatives act only on $\x$, the parameters behave as additional coordinates without derivatives, which is how they are discretized in Sec.~\ref{sec:qtt-parametric}.

The resulting variational structure makes the problem naturally suited to tensor networks. A conforming Galerkin discretization yields a symmetric positive-definite (SPD) linear system whose solution minimizes the discrete energy. In compressed form, this minimization is closely analogous to a ground-state search, with the discrete operator in the role of a quantum Hamiltonian. The variational sweeping algorithms developed for ground-state search, such as the density-matrix renormalization group (DMRG)~\cite{white1992,schollwock2011}, therefore carry over, with local eigenvalue problems replaced by local linear systems, as in the alternating linear scheme and its modified variant (ALS/MALS)~\cite{holtz2012} and the alternating minimal energy method (AMEn)~\cite{dolgov2014}. Since QTT truncation perturbs the exact operator, with effects amplified by poor conditioning, the truncation tolerance must be tight enough to preserve both the accuracy and the SPD structure. Nonsymmetric problems can be treated with TT-GMRES~\cite{dolgov2012b} or nonsymmetric alternating schemes~\cite{dolgov2014}, but without a convex energy they lose this minimization structure and may be less effective in practice~\cite{rohrig2025}; we therefore focus on the symmetric, coercive setting.

\section{QTT Chebyshev--Galerkin Method}
\label{sec:qtt}

\begin{figure}
    \centering
    \includegraphics[width=\linewidth]{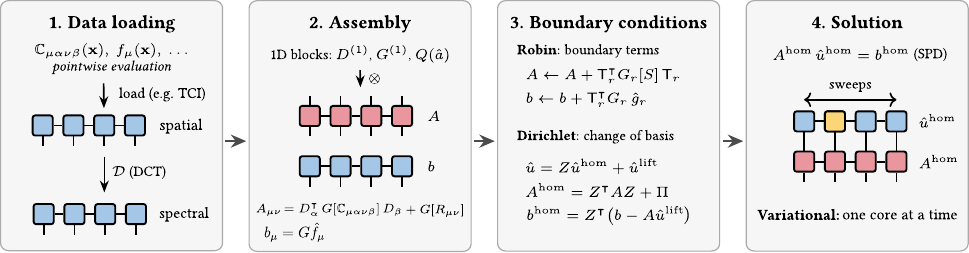}
    \caption{Overview of the QTT Chebyshev--Galerkin method. (1) Problem data are represented as nodal QTTs on the Chebyshev grid and mapped to Chebyshev coefficients by the transform MPO. (2) The Galerkin operator and load vector are assembled from one-dimensional differentiation, Gram, and multiplication MPOs. (3) Robin conditions add boundary terms, and Dirichlet conditions are imposed by basis recombination and lifting. (4) The resulting SPD system is solved by a variational QTT solver and the solution is evaluated from its Chebyshev coefficients. All stages operate in compressed QTT/MPO form.}
    \label{fig:pipeline}
\end{figure}

We now discretize the weak problem of Sec.~\ref{sec:problem} with a spectral Chebyshev--Galerkin method, performing every step in compressed form. The construction proceeds in four stages, summarized in Fig.~\ref{fig:pipeline}. First, the problem data, available only through pointwise evaluation, are loaded as QTTs on a grid, typically using tensor cross-interpolation (TCI)~\cite{oseledets2010cross,ritter2024}, and transformed into Chebyshev coefficients. Second, the Galerkin operator and load vector are assembled from a few one-dimensional building blocks---differentiation, integration, and multiplication by the problem data---each represented as an MPO. Third, Robin conditions are added as boundary terms, and Dirichlet conditions are imposed through a change of basis. Finally, the resulting SPD system is solved with a variational tensor-network solver. No stage forms the full tensor-product arrays, so the cost is governed throughout by bond dimensions. Section~\ref{sec:qtt-qtt} reviews the QTT format, Sec.~\ref{sec:qtt-galerkin} develops a basis-independent Galerkin construction, including boundary conditions and the parametric extension, Sec.~\ref{sec:qtt-chebyshev} specializes it to Chebyshev polynomials, and Sec.~\ref{sec:qtt-solve} addresses the solution of the resulting linear system.

\subsection{QTT and MPO Representations}
\label{sec:qtt-qtt}

\begin{figure}
    \centering
    \includegraphics[width=0.75\linewidth]{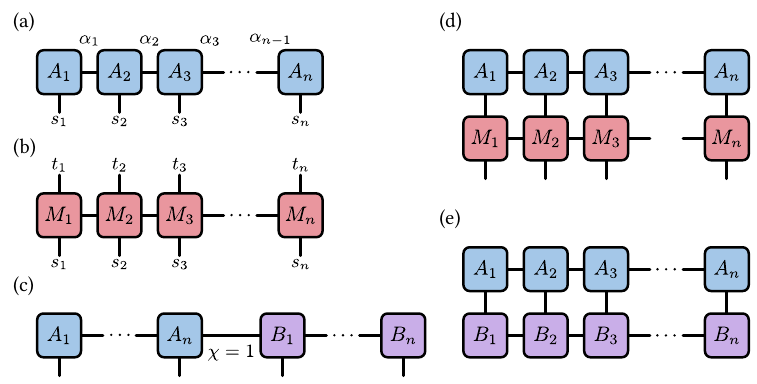}
    \caption{QTT and MPO representations and basic operations. (a) QTT with physical indices $s_k$ and virtual indices $\alpha_k$. (b) MPO with row indices $s_k$ and column indices $t_k$. (c) Kronecker product $A\otimes B$, obtained by concatenating two chains through a unit bond. (d) MPO--QTT product $MA$, obtained by contracting the shared physical indices. (e) Inner product $A\cdot B$.}
    \label{fig:tt}
\end{figure}

A tensor train (TT) factorizes an array with $n$ indices into a chain of $n$ small cores~\cite{oseledets2011},
\begin{equation}
    A_{s_1 s_2 \ldots s_n}
    = \sum_{\alpha_1,\ldots,\alpha_{n-1}}
    A^{(1)}_{s_1\alpha_1}\, A^{(2)}_{\alpha_1 s_2 \alpha_2} \cdots A^{(n)}_{\alpha_{n-1} s_n},
    \label{eq:tt}
\end{equation}
as depicted in Fig.~\ref{fig:tt}(a). The physical indices $s_k$ are those of the original array, while the virtual indices $\alpha_k=1,\ldots,\chi_k$ connect neighboring cores. The bond dimensions $\chi_k$ control both accuracy and cost: with physical dimension $\ell$ and maximum bond dimension $\chi$, a TT stores $\O{n\ell\chi^2}$ numbers instead of the $\ell^n$ entries of the full array.

In the quantics tensor train (QTT) format~\cite{oseledets2010,khoromskij2011b}, the physical indices are the binary digits of a grid index. For a function $a(x)$ sampled at $N=2^n$ grid points $x_j$, the grid index is written in binary,
\begin{equation}
    j = \sum_{k=1}^{n} s_k\, 2^{n-k}, \qquad s_k\in\{0,1\},
    \label{eq:qtt_binary}
\end{equation}
so that $a(x_j) = a_{s_1 s_2 \ldots s_n}$, where $s_1$ is the most significant bit. The resulting array with $n$ binary indices is then represented as a TT with $n$ cores of physical dimension two. For bounded bond dimension, storage therefore grows linearly with $n=\log_2 N$, that is, logarithmically with resolution. In several dimensions, the bits of each coordinate are placed contiguously along the chain~\cite{garcia2021}, so that each coordinate contributes its own block of cores; further indices, such as field components or parameters, simply append further cores. How well a function compresses depends on its smoothness and structure, as well as on the ordering of the bits~\cite{shinaoka2023,aldavero2024}.

Linear operators are represented analogously as matrix-product operators (MPOs), whose cores carry both a row and a column index [Fig.~\ref{fig:tt}(b)]. QTTs and MPOs support the operations of linear algebra directly on their cores~\cite{garcia2021,garcia2026}: sums, elementwise products, matrix--vector and matrix--matrix products [Fig.~\ref{fig:tt}(d)], and inner products and integrals [Fig.~\ref{fig:tt}(e)]. Kronecker products are particularly simple: they amount to concatenating two chains through a bond of dimension one [Fig.~\ref{fig:tt}(c)], which involves no arithmetic and no increase in bond dimension, whereas the explicit Kronecker product of two arrays of sizes $N_1$ and $N_2$ requires $\mathcal{O}(N_1N_2)$ operations. Multidimensional operators built from one-dimensional QTT blocks are therefore essentially free to assemble. The remaining operations generally increase the bond dimension, which is then reduced by truncation: each bond is compressed through a singular value decomposition, discarding singular values below a relative tolerance $\varepsilon$~\cite{oseledets2011}. QTT arithmetic is therefore a form of finite-precision linear algebra, in which $\varepsilon$ trades accuracy for compression. Unless stated otherwise, all vectors and operators below are QTTs and MPOs, respectively.

\subsection{Spectral Galerkin Discretization}
\label{sec:qtt-galerkin}

The Galerkin method replaces the infinite-dimensional spaces of Eq.~\eqref{eq:spaces} by the span of a finite basis and requires the weak form, Eq.~\eqref{eq:bvp_weak}, to hold for every basis function. We first develop the construction for a generic one-dimensional basis, such as Legendre polynomials~\cite{shen1994} or Fourier modes~\cite{shen2011}, and specialize it to the Chebyshev basis~\cite{shen1995} in Sec.~\ref{sec:qtt-chebyshev}. Dirichlet conditions are imposed at the end of this section, so trial and test functions are unconstrained until then.

\subsubsection{Tensor-product basis}
Along each coordinate $x_\gamma$, let $\{\phi_k\}_{k=0}^{N_\gamma-1}$ be a set of $N_\gamma$ one-dimensional basis functions. The $d$-dimensional basis consists of their products,
\begin{equation}
    \Phi_{\bs{k}}(\x) = \prod_{\gamma=1}^{d}\phi_{k_\gamma}(x_\gamma),
    \label{eq:basis_tensor}
\end{equation}
labeled by the multi-index $\bs{k}=(k_1,\ldots,k_d)$; boldface indices denote multi-indices throughout. The resulting space contains $N=\prod_\gamma N_\gamma$ basis functions, and each field component is expanded as
\begin{equation}
    u_\nu(\x) \approx \sum_{\bs{k}} \hat{u}_{\nu\bs{k}}\, \Phi_{\bs{k}}(\x).
    \label{eq:u_expansion}
\end{equation}
For each component, the spectral coefficients $\hat{u}_{\nu\bs{k}}$ form a vector $\hat{u}_\nu\in\mathbb{R}^N$, so that the complete unknown $\hat{u}=(\hat{u}_1,\ldots,\hat{u}_m)$ belongs to $\mathbb{R}^m\otimes\mathbb{R}^N$. In QTT form, each coordinate $x_\gamma$ contributes a block of $n_\gamma = \log_2 N_\gamma$ cores, and the component index $\nu$ one additional core.

\subsubsection{Galerkin system}
Taking as test functions $v_\mu=\Phi_{\bs{i}}$, for every component $\mu$ and basis function $\bs{i}$, and substituting Eq.~\eqref{eq:u_expansion} into Eq.~\eqref{eq:bvp_weak} yields the linear system
\begin{equation}
    \sum_{\nu=1}^m A_{\mu\nu}\,\hat{u}_\nu = b_\mu, \qquad \mu=1,\ldots,m,
    \label{eq:galerkin_system}
\end{equation}
whose blocks $A_{\mu\nu}\in\mathbb{R}^{N\times N}$ and load vectors $b_\mu\in\mathbb{R}^N$ have entries
\begin{equation}
  \begin{aligned}
    \big[A_{\mu\nu}\big]_{\bs{i}\bs{j}}
    &= \io \mathbb{C}_{\mu\alpha\nu\beta}(\x)\, \p_\alpha\Phi_{\bs{i}}(\x)\, \p_\beta\Phi_{\bs{j}}(\x)\, \d\x
    + \io R_{\mu\nu}(\x)\, \Phi_{\bs{i}}(\x)\Phi_{\bs{j}}(\x)\, \d\x
    + \int_{\p\Omega} S_{\mu\nu}(\x)\, \Phi_{\bs{i}}(\x)\Phi_{\bs{j}}(\x)\, \dS,
    \\
    \big[b_\mu\big]_{\bs{i}}
    &= \io f_\mu(\x)\, \Phi_{\bs{i}}(\x)\, \d\x
    + \int_{\p\Omega} g_\mu(\x)\, \Phi_{\bs{i}}(\x)\, \dS.
  \end{aligned}
  \label{eq:galerkin_discrete}
\end{equation}
Every term has the same structure: the integral of a problem coefficient times two basis functions, possibly differentiated. Assembling the system therefore requires only three operations on spectral coefficients---differentiating a function, multiplying it by a coefficient field, and integrating the product of two functions---which we now express as matrices.

\subsubsection{One-dimensional building blocks}
Consider a single coordinate, and let a function $v(x)=\sum_k \hat{v}_k\phi_k(x)$ be represented by its coefficient vector $\hat{v}$.
\begin{itemize}
  \item \emph{Differentiation.} The differentiation matrix $D^{(1)}$ maps the coefficients of $v$ to those of its derivative $v'$. It is defined by its action on the basis,
  \begin{equation}
      \phi_k'(x) = \sum_j D^{(1)}_{jk}\,\phi_j(x).
  \end{equation}
  For polynomial bases, differentiation lowers the degree, so the retained space is closed under differentiation; the unused highest mode is padded with zeros to keep $D^{(1)}$ square.

  \item \emph{Multiplication.} The multiplication matrix $Q(a)$ maps the coefficients of $v$ to those of the product $a(x)v(x)$, where $a(x)=\sum_k\hat{a}_k\phi_k(x)$ is a coefficient field of the problem. It is defined by
  \begin{equation}
      a(x)\,\phi_k(x) = \sum_j Q(a)_{jk}\,\phi_j(x),
      \label{eq:multiplication_operator}
  \end{equation}
  and depends on $a$ only through its coefficients $\hat{a}$, following the product rule of the basis functions (Sec.~\ref{sec:qtt-chebyshev} gives it for Chebyshev polynomials). Because the product $a\phi_k$ contains higher modes than either factor, $Q(a)$ is formed on an enlarged space with $2N_\gamma$ modes, padding coefficient vectors with zeros, so that the product is represented exactly before integration.

  \item \emph{Integration.} The Gram matrix, or overlap matrix,
  \begin{equation}
      G^{(1)}_{ij} = \int_{-1}^{1}\phi_i(x)\phi_j(x)\,\d x
  \end{equation}
  turns integrals of products into products of coefficient vectors, $\int v\,w\,\d x = \hat{v}^\intercal G^{(1)}\hat{w}$. For several polynomial bases, $G^{(1)}$ is not diagonal, since their orthogonality holds with respect to a weight function rather than the plain integral.
\end{itemize}

Combining multiplication and integration gives the weighted Gram matrix, which contains the integrals of the coefficient field against pairs of basis functions:
\begin{equation}
    G[a]_{ij} \coloneqq \int_{-1}^{1} a(x)\,\phi_i(x)\phi_j(x)\,\d x
    = \sum_{k} G^{(1)}_{ik}\, Q(a)_{kj},
    \label{eq:weighted_gram}
\end{equation}
with the sum over $k$ extending over the enlarged space, i.e., $G[a] = P^\intercal G^{(1)}Q(a)\,P$, where $P:\mathbb{R}^{N_\gamma}\to\mathbb{R}^{2N_\gamma}$ pads vectors with zeros to the enlarged space, on which $G^{(1)}$ and $Q(a)$ act. The weighted integral is thus computed by first multiplying by $a$ in coefficient space and then integrating against the test basis. This factorization isolates all dependence on the problem data in $Q(a)$, while $G^{(1)}$ depends only on the basis. Although written as a product, $G[a]$ is symmetric, as follows from its definition.

\subsubsection{Multidimensional operators}
For the tensor-product basis of Eq.~\eqref{eq:basis_tensor}, differentiation and integration factorize into one-dimensional blocks,
\begin{equation}
    D_\alpha = I_1\otimes\cdots\otimes I_{\alpha-1}\otimes D^{(1)}_\alpha \otimes I_{\alpha+1}\otimes\cdots\otimes I_d,
    \qquad
    G = \bigotimes_{\gamma=1}^d G^{(1)}_\gamma,
    \label{eq:DG_kronecker}
\end{equation}
where $D^{(1)}_\alpha$ and $G^{(1)}_\gamma$ are the one-dimensional matrices along $x_\alpha$ and $x_\gamma$, and $I_\gamma$ is the identity on the remaining coordinates: differentiation along $x_\alpha$ leaves all other coordinates untouched, whereas integration acts on all of them. Since Kronecker products are essentially free in QTT form (Sec.~\ref{sec:qtt-qtt}), both operators follow directly from their one-dimensional MPOs. The multiplication operator, in contrast, is generally not a Kronecker product, because coefficient fields such as a spatially varying elastic modulus need not be separable. It is defined as in Eq.~\eqref{eq:multiplication_operator}, with $\Phi_{\bs{k}}$ in place of $\phi_k$, and constructed directly from the QTT of the multivariate coefficients $\hat{a}_{\bs{k}}$ (see Sec.~\ref{sec:qtt-chebyshev}). The weighted Gram matrix then follows as
\begin{equation}
    G[a] = P^\intercal G\,Q(a)\,P, \qquad P=\bigotimes_\gamma P_\gamma,
    \label{eq:weighted_gram_decomposition}
\end{equation}
where $P_\gamma$ pads the coefficients with zeros along $x_\gamma$.

\subsubsection{Bulk terms}
With these operators, the bulk integrals in Eq.~\eqref{eq:galerkin_discrete} become
\begin{equation}
    A^{\mathrm{bulk}}_{\mu\nu} = \sum_{\alpha,\beta=1}^d D_\alpha^\intercal\, G\big[\mathbb{C}_{\mu\alpha\nu\beta}\big]\, D_\beta
    + G\big[R_{\mu\nu}\big],
    \qquad
    b^{\mathrm{bulk}}_\mu = G\,\hat{f}_\mu.
    \label{eq:A_bulk}
\end{equation}
Read from right to left, the first term of $A^{\mathrm{bulk}}_{\mu\nu}$ differentiates the trial function along $x_\beta$, multiplies it by the coefficient field $\mathbb{C}_{\mu\alpha\nu\beta}$, and integrates the result against the test function differentiated along $x_\alpha$; the second term is analogous but without derivatives. In the load vector, the source enters through its coefficient vector $\hat{f}_\mu$, since $\io f_\mu\Phi_{\bs{i}}\,\d\x = (G\hat{f}_\mu)_{\bs{i}}$.

\subsubsection{Robin boundary conditions}
\label{sec:qtt-robin}
The boundary integrals are evaluated face by face. Each Robin face $r=1,\ldots,n_\mathrm{R}$ lies at $x_{\gamma}=\sigma$, with $\sigma\in\{-1,+1\}$, and restricting a basis function to it only evaluates the factor along $x_\gamma$. The trace operator onto that face is therefore
\begin{equation}
    \mathsf{T}_r = I_1\otimes\cdots\otimes t_\sigma^\intercal\otimes\cdots\otimes I_d,
    \qquad
    t_\sigma = \big(\phi_0(\sigma),\ldots,\phi_{N_\gamma-1}(\sigma)\big)^\intercal,
    \label{eq:trace_operator}
\end{equation}
a rectangular matrix that maps the $N$ coefficients of a function on $\Omega$ to the $N/N_\gamma$ coefficients of its restriction to the face. In a nodal representation, $\mathsf{T}_r$ would simply select the grid points on the face; in coefficient space, it instead combines all modes along $x_\gamma$ through their boundary values $\phi_k(\sigma)$. The Robin contributions are then
\begin{equation}
    A^{\mathrm{Robin}}_{\mu\nu} = \sum_{r=1}^{n_\mathrm{R}} \mathsf{T}_r^\intercal\, G_r\big[S_{\mu\nu}\big]\, \mathsf{T}_r,
    \qquad
    b^{\mathrm{Robin}}_\mu = \sum_{r=1}^{n_\mathrm{R}} \mathsf{T}_r^\intercal\, G_r\, \hat{g}_{\mu,r},
    \label{eq:A_robin}
\end{equation}
where $G_r[a]$ is the weighted Gram matrix on the face and $\hat{g}_{\mu,r}$ contains the coefficients of the Robin data there; both vanish on faces where component $\mu$ is not of Robin type. Each term restricts trial and test functions to the face, integrates them there, and maps the result back to the full coefficient space. Robin conditions therefore enter as additional terms built from the same operations as the bulk, without modifying any entry of the bulk operator.

Adding both contributions gives the blocks of Eq.~\eqref{eq:galerkin_system},
\begin{equation}
    A_{\mu\nu} = A^{\mathrm{bulk}}_{\mu\nu} + A^{\mathrm{Robin}}_{\mu\nu},
    \qquad
    b_\mu = b^{\mathrm{bulk}}_\mu + b^{\mathrm{Robin}}_\mu,
    \label{eq:A_galerkin}
\end{equation}
which are combined into a single operator and load vector,
\begin{equation}
    A = \sum_{\mu,\nu=1}^m E_{\mu\nu}\otimes A_{\mu\nu},
    \qquad
    b = \sum_{\mu=1}^m \mathbf{e}_\mu\otimes b_\mu,
    \qquad
    E_{\mu\nu}=\mathbf{e}_\mu\mathbf{e}_\nu^\intercal,
    \label{eq:vector_assembly}
\end{equation}
where $\mathbf{e}_\mu$ is the $\mu$-th canonical basis vector of $\mathbb{R}^m$, so that Eq.~\eqref{eq:galerkin_system} reads $A\hat{u}=b$ on $\mathbb{R}^m\otimes\mathbb{R}^N$. In QTT form, the factors $E_{\mu\nu}$ and $\mathbf{e}_\mu$ form the additional core that carries the component index. Under the symmetry assumptions of Eq.~\eqref{eq:bvp_symmetry}, $A_{\mu\nu}=A_{\nu\mu}^\intercal$, so that $A$ is symmetric.

\subsubsection{Dirichlet boundary conditions}
\label{sec:qtt-dirichlet}

Dirichlet conditions are imposed through the trial and test spaces of Eq.~\eqref{eq:spaces} rather than through additional terms. Following Sec.~\ref{sec:problem-weak}, the solution is split into a lifting function, which satisfies the prescribed data, and a homogeneous correction, which vanishes on the Dirichlet faces:
\begin{equation}
    \hat{u} = \ulift + Z\,\uhom.
    \label{eq:u_dirichlet_decomposition}
\end{equation}
Here, $\ulift$ contains the coefficients of any function that matches $\bar{u}_\mu$ on $\GDirichlet_\mu$, and $Z$ is a recombination matrix whose columns are combinations of basis functions that vanish on the Dirichlet faces, so that every vector $\uhom$ yields an admissible correction.

Along a single coordinate, a combination of basis functions with coefficient vector $z$ vanishes at $x=\sigma$ precisely when $t_\sigma^\intercal z = 0$, with $t_\sigma$ the vector of boundary values in Eq.~\eqref{eq:trace_operator}. The one-dimensional recombination matrix $Z^{(1)}$ collects $N_\gamma-1$ such vectors for a constraint at one endpoint, or $N_\gamma-2$ for constraints at both. For Chebyshev polynomials, which satisfy $T_k(\pm1)=(\pm1)^k$, simple choices are
\begin{equation}
    \phi^{\sigma}_k(x) = T_{k+1}(x) - \sigma\,T_k(x),
    \qquad
    \phi^{0}_k(x) = T_{k+2}(x) - T_k(x),
\end{equation}
which vanish at $x=\sigma$ and at both endpoints, respectively. Their coefficient vectors form the columns of
\begin{equation}
    Z^{(1)}_{jk} =
    \begin{cases}
        \delta_{j\,k+1} - \sigma\,\delta_{jk}, & \text{constraint at } x=\sigma, \quad k=0,\ldots,N_\gamma-2,
        \\[0.4em]
        \delta_{j\,k+2} - \delta_{jk}, & \text{constraints at both endpoints,} \quad k=0,\ldots,N_\gamma-3,
        \\[0.4em]
        0, & \text{remaining columns},
    \end{cases}
    \label{eq:recombination_matrix}
\end{equation}
with $Z^{(1)}=I$ along unconstrained coordinates. Each nonzero column has only two nonzero entries, so $Z^{(1)}$ is a banded matrix with a compact MPO representation. The zero columns pad $Z^{(1)}$ to size $N_\gamma\times N_\gamma$, which preserves power-of-two mode sizes in QTT form. For component $\mu$, the multidimensional recombination matrix is the Kronecker product $Z_\mu=\bigotimes_\gamma Z^{(1)}_{\mu\gamma}$, and $Z=\sum_\mu E_{\mu\mu}\otimes Z_\mu$ acts on all components.

Substituting Eq.~\eqref{eq:u_dirichlet_decomposition} into $A\hat{u}=b$ and testing only with recombined functions gives the reduced system
\begin{equation}
    \Ahom\,\uhom = \bhom,
    \qquad
    \Ahom = Z^\intercal A Z + \Pi,
    \qquad
    \bhom = Z^\intercal\big(b - A\ulift\big),
    \label{eq:galerkin_system_lifted}
\end{equation}
where the diagonal matrix $\Pi$ acts only on the padded coefficients and fixes them to zero, with entries scaled consistently with $Z^\intercal AZ$ to avoid degrading the conditioning. Since $Z^\intercal AZ$ represents the bilinear form on the discrete test space, $\Ahom$ is SPD under the coercivity assumption of Sec.~\ref{sec:problem-weak}. Homogeneous conditions correspond to $\ulift=0$, while nonhomogeneous lifts are constructed directly in coefficient space by spectral transfinite interpolation (see Appendix~\ref{app:tfi}). As with Robin conditions, Dirichlet conditions are thus imposed through operations on whole operators---here a change of basis---rather than by modifying individual entries.

\subsubsection{Spatial--parametric extension}
\label{sec:qtt-parametric}

The parameters $\bs{\xi}\in\Xi$ of Eq.~\eqref{eq:bvp_weak_parametric} are discretized exactly like the spatial coordinates. Along each parameter $\xi_q$, we introduce $N_\xi$ one-dimensional basis functions $\psi_\ell$ and expand the solution in the joint tensor-product basis,
\begin{equation}
    u_\nu(\x,\bs{\xi}) \approx \sum_{\bs{k},\bs{\ell}} \hat{u}_{\nu\bs{k}\bs{\ell}}\, \Phi_{\bs{k}}(\x)\,\Psi_{\bs{\ell}}(\bs{\xi}),
    \qquad
    \Psi_{\bs{\ell}}(\bs{\xi}) = \prod_{q=1}^{d_\xi}\psi_{\ell_q}(\xi_q).
    \label{eq:u_expansion_parametric}
\end{equation}
Testing Eq.~\eqref{eq:bvp_weak_parametric} with the products $\Phi_{\bs{i}}\Psi_{\bs{\ell}'}$ yields a linear system with the same block structure as Eq.~\eqref{eq:galerkin_system}, now on $\mathbb{R}^m\otimes\mathbb{R}^N\otimes\mathbb{R}^{N_\Xi}$ with $N_\Xi=N_\xi^{d_\xi}$. Every ingredient of the construction carries over to the enlarged basis. The Gram matrix acquires one factor per parameter, $G\to G\otimes G_\Xi$ with $G_\Xi=\bigotimes_q G^{(1)}_{\xi_q}$, which realizes the integral over $\Xi$ in Eq.~\eqref{eq:bvp_weak_parametric}. The multiplication operator is built from the coefficients of $a(\x,\bs{\xi})$ over all coordinates, and thereby captures the joint dependence of the data on space and parameters. Traces and recombination act on the spatial factors only, $\mathsf{T}_r\to\mathsf{T}_r\otimes I_\Xi$ and $Z\to Z\otimes I_\Xi$, while the lift may depend on $\bs{\xi}$ through the Dirichlet data. The only essential difference is that the weak form contains no derivatives with respect to $\bs{\xi}$, so each differentiation matrix is extended by the identity, $D_\alpha \to D_\alpha\otimes I_\Xi$.

Parameters therefore behave as additional coordinates along which every differentiation operator acts as the identity. In QTT form, each parameter contributes a further block of $n_\xi = \log_2 N_\xi$ cores, and a single compressed representation encodes the complete family $u(\x,\bs{\xi})$. One solve thus returns the solution for every parameter value, rather than one solve per sample.

\subsection{Chebyshev--Galerkin Discretization}
\label{sec:qtt-chebyshev}

We now specialize the construction of Sec.~\ref{sec:qtt-galerkin} to the basis of Chebyshev polynomials of the first kind. Three properties make it particularly suited to QTT form: Chebyshev expansions converge rapidly for smooth functions~\cite{trefethen2000,boyd2001}, point values and coefficients are related by a fast discrete cosine transform~\cite{boyd2001}, and the algebra of the polynomials yields closed-form, highly structured differentiation, Gram, and multiplication matrices, all of which admit compact MPO representations.

\subsubsection{Tensor-product Chebyshev basis}
\begin{figure}
    \centering
    \includegraphics[width=0.8\linewidth]{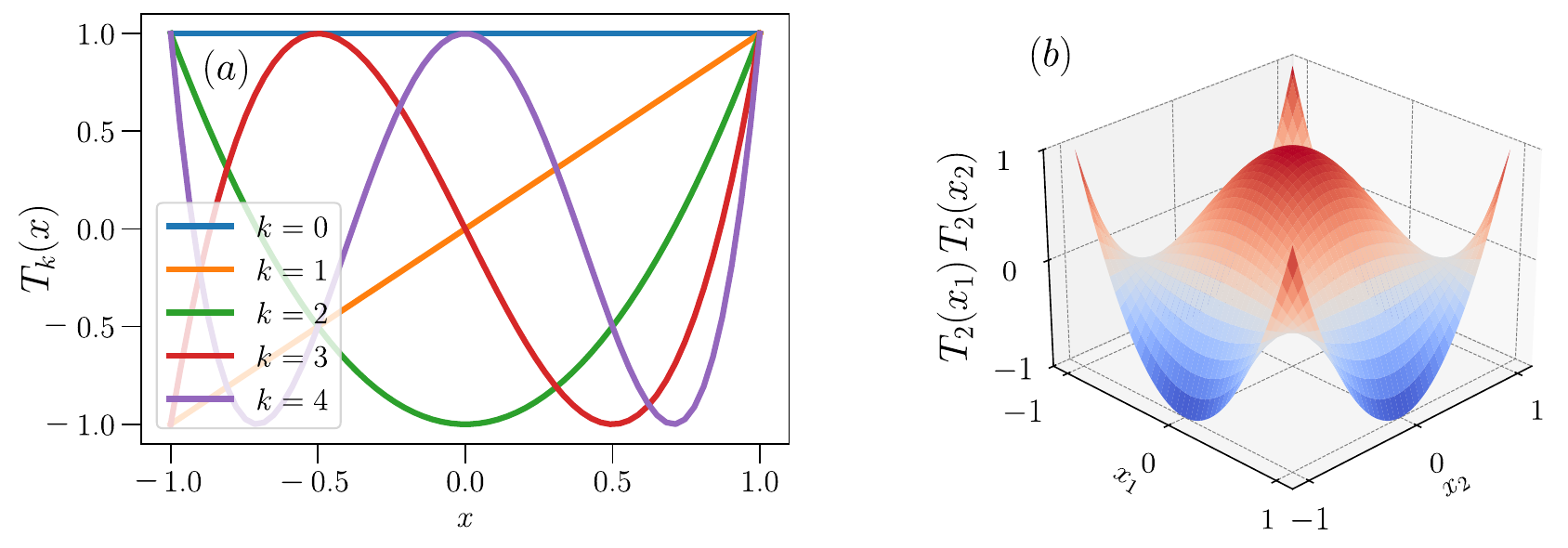}
    \caption{(a) First five Chebyshev polynomials of the first kind, $T_k(x)$, for $k=0,\ldots,4$ on $x\in[-1,1]$. (b) Tensor-product basis function $T_2(x_1)T_2(x_2)$ on $[-1,1]^2$.}
    \label{fig:chebyshev}
\end{figure}

On the reference interval $[-1,1]$, the Chebyshev polynomials of the first kind are
\begin{equation}
    T_k(x) = \cos(k\arccos x), \qquad k\geq 0,
\end{equation}
shown in Fig.~\ref{fig:chebyshev}. We take $\phi_k=T_k$ in Eq.~\eqref{eq:basis_tensor}, so that the multidimensional basis functions are
\begin{equation}
    \T_{\bs{k}}(\x) = \prod_{\gamma=1}^{d}T_{k_\gamma}(x_\gamma),
    \label{eq:chebyshev_tensor}
\end{equation}
and retain $N_\gamma=2^{n_\gamma}$ modes per coordinate, so that each modal index $k_\gamma$ is represented by $n_\gamma$ binary digits. Since $T_k(\pm1)=(\pm1)^k$, the boundary values entering the trace operator of Eq.~\eqref{eq:trace_operator} are simply $t_\sigma=(1,\sigma,\sigma^2,\ldots)^\intercal$, $\sigma = \pm 1$.

\subsubsection{Loading the problem data into Chebyshev space}
The construction requires the Chebyshev coefficients $\hat{a}_{\bs{k}}$ of every problem input---the coefficient fields, sources, and boundary data---in QTT form. Since these inputs are typically given as functions that can be evaluated pointwise, their coefficients are obtained in two steps. 

First, each input is represented as a nodal QTT of its values on the tensor-product Chebyshev--Gauss grid, with nodes
\begin{equation}
    x_j = \cos\theta_j, \qquad \theta_j = \frac{(2j+1)\pi}{2N_\gamma}, \qquad j=0,\ldots,N_\gamma-1,
\end{equation}
along each coordinate. These nodes, like the polynomials themselves, are defined on $[-1,1]$, which is why the problem is posed on the reference domain $\Omega=[-1,1]^d$. Simple functions admit exact QTT constructions, whereas general data can be loaded by TCI~\cite{oseledets2010cross,ritter2024}, which builds the QTT from a small, adaptively chosen set of function evaluations without forming the full grid. 

Second, the nodal values are mapped to Chebyshev coefficients. Since $T_k(x_j)=\cos(k\theta_j)$ and the Chebyshev polynomials are discretely orthogonal on these nodes, the interpolating expansion $a(x)\approx\sum_k \hat{a}_k T_k(x)$ has coefficients
\begin{equation}
    \hat{a}_k = \frac{2-\delta_{k0}}{N_\gamma} \sum_{j=0}^{N_\gamma-1} a(x_j)\, \cos\!\left(\frac{k(2j+1)\pi}{2N_\gamma}\right),
    \qquad k=0,\ldots,N_\gamma-1.
    \label{eq:chebyshev_transform}
\end{equation}
Up to the diagonal factor $(2-\delta_{k0})/N_\gamma$, this is the discrete cosine transform of type II (DCT-II), the cosine counterpart of the discrete Fourier transform. We denote by $\mathcal{D}_\gamma$ the matrix of Eq.~\eqref{eq:chebyshev_transform} along coordinate $x_\gamma$, so that the coefficients of a multivariate input follow from its nodal QTT $a$ as
\begin{equation}
    \hat{a} = \mathcal{D}\,a,
    \qquad
    \mathcal{D}=\bigotimes_{\gamma=1}^d\mathcal{D}_\gamma.
\end{equation}
Like the Fourier transform, the Chebyshev transform in Eq.~\eqref{eq:chebyshev_transform} factors into $\mathcal{O}(\log_2 N_\gamma)$ steps, each acting on a single binary digit of the index~\cite{perera2016}. In QTT form, each step is an exact MPO with bond dimension at most 10, independently of the resolution (see Appendix~\ref{app:transform}).

\subsubsection{One-dimensional Chebyshev blocks}
With $\phi_k=T_k$, the three one-dimensional blocks of Sec.~\ref{sec:qtt-galerkin} take closed forms, each with a compact MPO representation (see Appendix~\ref{app:mpos}).
\begin{itemize}
  \item \emph{Differentiation matrix.} The derivative of a Chebyshev expansion has coefficients
  \begin{equation}
      D^{(1)}_{ij} =
      \begin{cases}
          j, & i=0 \text{ and } j \text{ odd},
          \\[0.4em]
          2j, & 1\leq i<j \text{ and } j-i \text{ odd},
          \\[0.4em]
          0, & \text{otherwise},
      \end{cases}
  \end{equation}
  an upper-triangular parity pattern that is represented exactly by an MPO of bond dimension 6 (see Appendix~\ref{app:mpos-differentiation}). Since its entries grow with the degree, the associated differential operator becomes ill-conditioned at high resolution, which we mitigate by preconditioning (see Appendix~\ref{app:preconditioning}).

  \item \emph{Gram matrix.} Integrating products of Chebyshev polynomials gives
  \begin{equation}
      G^{(1)}_{ij} =
      \begin{cases}
          \dfrac{1}{1-(i+j)^2} + \dfrac{1}{1-(i-j)^2}, & i+j \text{ even},
          \\[1em]
          0, & i+j \text{ odd}.
      \end{cases}
      \label{eq:gram_chebyshev}
  \end{equation}
  Each entry depends only on the sum and the difference of its indices, so $G^{(1)}$ is the sum of a Hankel and a Toeplitz matrix, both generated by the single vector $w_k=2/(1-k^2)$ for even $k$ and $w_k=0$ for odd $k$. This structure allows the MPO to be built directly from a QTT of $w$, without forming or contracting any matrix: each core of the QTT is mapped to an MPO core, augmented with a few auxiliary bond states that carry the sum $i+j$, or the comparison underlying $i-j$, from one binary digit to the next. The resulting MPO is exact for the given QTT, with bond dimension $\chi_G\leq6\chi_w+1$. The Gram matrix is thus constructed without quadrature: its only approximation is the TT-SVD compression of the vector $w$, of length $2N_\gamma$ (see Appendix~\ref{app:mpos-gram}).

  \item \emph{Multiplication matrix.} The Chebyshev product rule $T_iT_j=\tfrac12\left(T_{i+j}+T_{|i-j|}\right)$ turns Eq.~\eqref{eq:multiplication_operator} into
  \begin{equation}
      Q(a)_{ij} = \frac12 \sum_{k=0}^{M_\gamma-1} \hat{a}_k
      \left(\delta_{i\,j+k} + \delta_{i\,|k-j|}\right),
      \qquad
      0\leq i,j<M_\gamma,
      \label{eq:Q_chebyshev}
  \end{equation}
  on a modal space of size $M_\gamma=2N_\gamma$, with $\hat{a}_k=0$ for $k\geq N_\gamma$. Each entry is again a coefficient of $a$ selected by the sum or the difference of the indices, so $Q(a)$ has the same Hankel and Toeplitz structure as the Gram matrix in Eq.~\eqref{eq:gram_chebyshev}, now generated by $\hat{a}$. Its MPO is therefore built directly from the QTT of $\hat{a}$, with bond dimension $\chi_Q\leq7\chi_{\hat{a}}$ (see Appendix~\ref{app:mpos-multiplication}). For multivariate coefficients, the construction acts coordinate by coordinate on the QTT of $\hat{a}_{\bs{k}}$, so that $Q(a)$ inherits the possibly nonseparable structure of the data.
\end{itemize}

With these three blocks and the transform $\mathcal{D}$, every ingredient of Eqs.~\eqref{eq:A_galerkin}--\eqref{eq:vector_assembly} is available as an MPO: the basis-dependent blocks are exact or accurate to machine precision, and the problem data enter only through their coefficient QTTs $\hat{a}$. The complete Galerkin system is then assembled through Kronecker products, sums, and MPO products, without forming any tensor-product array, so that approximation errors arise only from loading the data and from truncating intermediate bond dimensions.

\subsection{Variational Solution}
\label{sec:qtt-solve}

Since the reduced system of Eq.~\eqref{eq:galerkin_system_lifted} is SPD, its solution is the unique minimizer of the discrete energy,
\begin{equation}
    \uhom = \operatorname*{arg\,min}_{v}\left\{\tfrac12\langle v,\Ahom v\rangle - \langle\bhom,v\rangle\right\}.
    \label{eq:discrete_energy}
\end{equation}
This problem closely resembles ground-state search in quantum many-body physics, where the Rayleigh quotient $\langle v,Hv\rangle/\langle v,v\rangle$ of a Hamiltonian $H$ is minimized over tensor networks. However, the two problems are not identical since the source term $\langle\bhom,v\rangle$ makes Eq.~\eqref{eq:discrete_energy} a linear problem rather than an eigenvalue problem. Still, the tensor-network machinery developed for computational quantum physics carries over naturally.

Variational solvers minimize Eq.~\eqref{eq:discrete_energy} directly in QTT form by sweeping back and forth along the chain and updating one core, or a pair of neighboring cores, while keeping the rest fixed. With the other cores fixed, the energy is quadratic in the free core, so each update reduces to a small linear system. Its size depends on the bond dimensions of $\Ahom$ and $\uhom$, not on the number of unknowns, and this is what keeps the cost polylogarithmic in the grid size. One option is MALS~\cite{holtz2012}, the inhomogeneous counterpart of the well-known DMRG algorithm~\cite{white1992}, which optimizes pairs of neighboring cores and adapts the bond dimension through SVD truncation. Instead, we use AMEn~\cite{dolgov2014}, which updates one core at a time and enriches it with a compressed approximation of the residual. Because its local problems involve a single core, AMEn is cheaper per sweep than two-core schemes~\cite{rohrig2025}, and the residual enrichment keeps it from stagnating on the ill-conditioned systems produced by the spectral discretization.

The spectral stiffness operator is ill-conditioned, and the local systems inherit this. We therefore precondition them with Jacobi scaling and, for the parametric problems, combine it with a Cholesky-based change of basis (Appendix~\ref{app:preconditioning}). Sweeps stop when the relative residual
\begin{equation}
    r = \frac{\|\Ahom\uhom-\bhom\|}{\|\bhom\|}
    \label{eq:residual}
\end{equation}
falls below a prescribed tolerance. This residual measures only the accuracy of the algebraic solve. The error from truncating the operator and the data is a separate contribution, which we assess against reference solutions in Sec.~\ref{sec:results}.

Once converged, the coefficients $\hat{u}=Z\uhom+\ulift$ define a global Chebyshev interpolant rather than values at a fixed set of nodes. A single QTT therefore encodes the field on all of $\Omega$. Its values on an equispaced grid of any density follow from applying the Chebyshev--Vandermonde matrix $V_{jk}=\T_k(y_j)$ as an MPO, with no loss of accuracy beyond that of the coefficients.

\subsection{Computational Complexity}
\label{sec:qtt-cost}

\begin{table}
    \centering
    \caption{Cost of the main operations in QTT form, with $n$ binary cores and bond dimensions $\chi_A$ and $\chi_u$, and for a dense spectral discretization with $\mathcal{N}=m\,2^n$ unknowns.}
    \label{tab:cost}
    \begin{tabular}{lll}
        \toprule
        Operation & QTT & Dense \\
        \midrule
        Storage of the solution & $\O{n\chi_u^2}$ & $\O{\mathcal{N}}$ \\
        Storage of the operator & $\O{n\chi_A^2}$ & $\O{\mathcal{N}^2}$ \\
        Operator--vector product & $\O{n\chi_A^3\chi_u^3}$ & $\O{\mathcal{N}^2}$ \\
        Chebyshev transform & $\O{n^2\chi_u^3}$ & $\O{\mathcal{N}\log\mathcal{N}}$ \\
        AMEn sweep & $\O{n\chi_A\chi_u^3 + n\chi_A^2\chi_u^2}$ & -- \\
        Linear solve & -- & $\O{\mathcal{N}^3}$ \\
        \bottomrule
    \end{tabular}
\end{table}

Let $n=\sum_\gamma\log_2 N_\gamma$ be the number of QTT cores, which grows logarithmically with the resolution and linearly with the number of spatial and parametric coordinates, and let $\mathcal{N}=m\,2^n$ be the nominal number of unknowns. Let $\chi_A$ and $\chi_u$ be the maximum bond dimensions of the operator and the solution. Table~\ref{tab:cost} compares the cost of the main operations in QTT form with that of a dense spectral discretization. In QTT form, every operation is linear in $n$, except the transform, whose $\O{\log_2 N_\gamma}$ layers per coordinate add a second factor of $n$ (see Appendix~\ref{app:transform}). For bounded bond dimensions, the cost of the whole method is therefore polylogarithmic in $\mathcal{N}$. In contrast, any conventional method that stores the unknowns explicitly requires at least $\O{\mathcal{N}}$ memory and time, merely to represent the solution, and a dense spectral solve requires $\O{\mathcal{N}^2}$ memory and $\O{\mathcal{N}^3}$ time. The QTT method thus has an exponential advantage in both time and memory over any explicit discretization with the same number of unknowns, provided that the bond dimensions remain bounded as the resolution increases. This condition is not guaranteed a priori, since the bond dimensions depend on the problem and may grow with the resolution. Moreover, the total cost depends on the number of solver sweeps and on the truncation tolerance needed for a given accuracy, both of which are affected by the conditioning of the problem. The experiments of Sec.~\ref{sec:results} therefore report the bond dimensions alongside the timings.

Accuracy, in turn, improves very rapidly with $n$. If the solution is analytic inside a Bernstein ellipse $E_\rho$ for some $\rho>1$, its Chebyshev coefficients decay as $\O{\rho^{-k}}$, and the error with $N_\gamma$ modes decreases as $\O{\rho^{-N_\gamma}}$~\cite{trefethen2000,boyd2001}. Since $N_\gamma=2^{n_\gamma}$, this discretization error then decays as $\O{\rho^{-2^{n_\gamma}}}$, doubly exponentially in the number of cores, so that each additional core roughly squares it. The computed QTT solution also carries truncation and solver errors, which set a floor on the attainable accuracy.

\section{Numerical Results}
\label{sec:results}

We test the method on four problems in linear elasticity, with the isotropic elasticity tensor of Eq.~\eqref{eq:elasticity_principal_tensor}. Each experiment supports one claim:
\begin{enumerate}
  \item \textbf{Manufactured verification (Sec.~\ref{sec:results-exp1}).} On a two-dimensional problem that combines a curved geometry, heterogeneous and coupled coefficients, and mixed boundary conditions, the QTT solution reproduces the spectral convergence of the dense Chebyshev--Galerkin discretization, down to an error floor set by the truncation tolerance.
  \item \textbf{High-resolution three-dimensional elasticity (Sec.~\ref{sec:results-exp2}).} On a three-dimensional problem with an oscillatory load, the method reaches errors below $10^{-8}$, which QTT finite differences and Q2/Q4 finite elements do not attain over the tested range, although QTT finite differences are cheaper at moderate accuracy.
  \item \textbf{High-dimensional parametric design (Sec.~\ref{sec:results-exp3}).} For a cantilever with six geometric, material, and boundary parameters, a single solve returns the response of $1.7\times10^{7}$ designs, with a compression factor of $6.5\times10^{4}$ over the explicit array.
  \item \textbf{Intrusive uncertainty quantification (Sec.~\ref{sec:results-exp4}).} For a block with a random elastic modulus described by ten random variables, a single solve returns the full random displacement field, with a compression factor of $2.4\times10^{6}$.
\end{enumerate}

All experiments can be reproduced with our implementation~\cite{rodriguez2026code}. Problem data are loaded by TCI~\cite{ritter2024,fernandez2024}, unless an exact QTT construction is available. Truncation is controlled by a relative tolerance $\varepsilon$ on the singular values at each bond, with no limit on the bond dimension. Since $\varepsilon$ does not bound the global error, all errors are measured against reference solutions. The linear systems are solved with AMEn, with the preconditioning described in Appendix~\ref{app:preconditioning}. Tolerances are given for each experiment. All computations were run on a workstation with an AMD Ryzen~9 5900X CPU (12 cores) and 64~GB of RAM.

\FloatBarrier
\subsection{Two-Dimensional Manufactured Verification}
\label{sec:results-exp1}

\begin{figure}
    \centering
    \includegraphics[width=\linewidth]{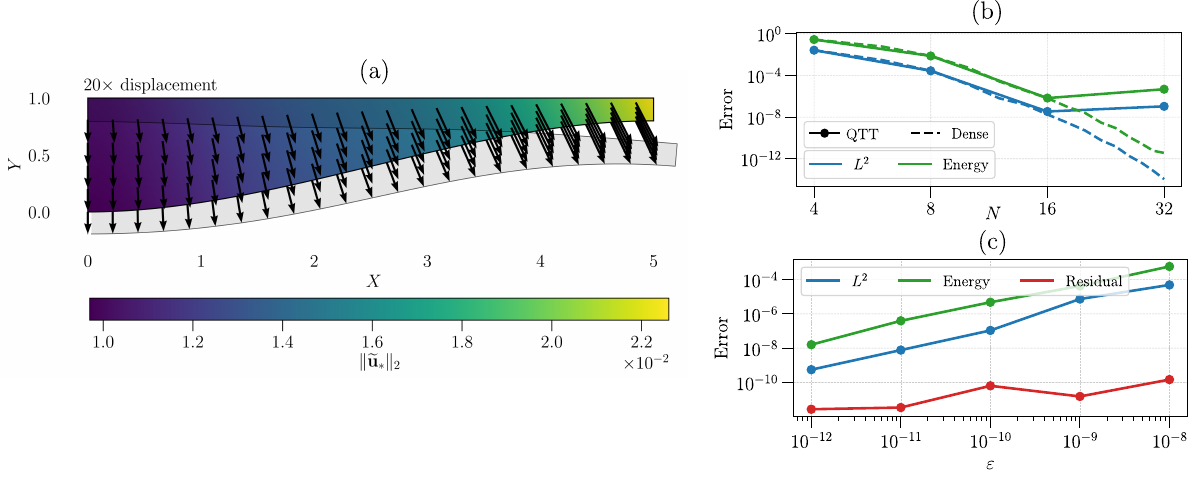}
    \caption{Manufactured verification. (a) Manufactured displacement on the tapered cantilever; color shows its magnitude, and vectors are magnified by a factor of 20. (b) Relative $\mathrm{L}^2$ and energy-norm errors of the dense and QTT discretizations versus $N$, with $\varepsilon=10^{-10}$. (c) Errors and AMEn residual versus truncation tolerance $\varepsilon$ at $N=32$. The point $\varepsilon=10^{-10}$ in (c) corresponds to the QTT result at $N=32$ in (b).}
    \label{fig:exp1}
\end{figure}

We first verify that the QTT formulation preserves the accuracy of the underlying Chebyshev--Galerkin discretization when all ingredients of the method act together. The problem is deliberately arbitrary but complete. It combines a curved geometry, a heterogeneous elastic modulus, coupled reaction terms in the bulk and on the boundary, and all three types of boundary conditions.

The domain is a tapered cantilever of length $\ell=5$, whose height decreases smoothly from 1 at the clamped end to $0.2$ at the tip [see Fig.~\ref{fig:exp1}(a)]. It is obtained from a smooth mapping of $\Omega=[-1,1]^2$ through
\begin{equation}
    X = \ell s,
    \qquad
    Y = t + (1-t)\,h(s),
    \qquad
    h(s) = 0.8\,s^2(3-2s),
\end{equation}
with $s=(1+x_1)/2$ and $t=(1+x_2)/2$, and pulled back as described in Appendix~\ref{app:mappings}. The material has a variable elastic modulus, $E(x_2)=1+t$, and a constant $\nu_{\mathrm P}=0.3$. The bulk and Robin terms couple both displacement components through
\begin{equation}
    R=(1+s)R_0,
    \qquad
    S=(1+t)R_0,
    \qquad
    R_0=\begin{pmatrix}1 & 1/2 \\ 1/2 & 1\end{pmatrix}.
    \label{eq:exp1-reactions}
\end{equation}
The manufactured displacement $\bs{u}_*$ combines axial stretching, bending, and Poisson contraction, so that all strain components are nonzero. Its bending profile is a static cantilever deflection perturbed by a higher vibration mode, which makes it smooth but not polynomial (see Appendix~\ref{app:manufactured}). The source and boundary data are computed from $\bs{u}_*$, with Dirichlet data on the clamped end, tractions on the upper and lower boundaries, and the Robin condition on the tip.

We compare the QTT solution for $N=2^n$, $n=2,\ldots,5$, modes per direction with a dense solution of the same discretization. We measure the relative $\mathrm{L}^2$ error on the Chebyshev--Gauss grid and the relative energy-norm error
\begin{equation}
    \mathrm{Err}_e = \left(\frac{\io \mathbb{C}_{\mu\alpha\nu\beta}\,\p_\alpha e_\mu\,\p_\beta e_\nu \,\d\x}{\io \mathbb{C}_{\mu\alpha\nu\beta}\,\p_\alpha u_{*\mu}\,\p_\beta u_{*\nu}\,\d\x}\right)^{1/2},
    \label{eq:elastic_error}
\end{equation}
where $e = u - u_*$, which measures the error in the elastic energy and is therefore sensitive to errors in strains and stresses.

The results in Fig.~\ref{fig:exp1} show that the QTT solution converges spectrally until truncation dominates, and that this truncation, not the solver, sets the final error. With $\varepsilon=10^{-10}$, the QTT and dense errors coincide up to $N=16$, where both are $\sim\!10^{-8}$ [see Fig.~\ref{fig:exp1}(b)]. At $N=32$, the dense error drops to $10^{-14}$, whereas the QTT error levels off at $1.1\times10^{-7}$, while the maximum bond dimension of the solution only grows from 4 at $N=4$ to 19 at $N=32$. At this resolution, the errors decrease steadily as $\varepsilon$ is tightened, while the AMEn residual of Eq.~\eqref{eq:residual} remains at or below about $10^{-10}$ for every tolerance [see Fig.~\ref{fig:exp1}(c)]. A small residual and a larger error are therefore consistent: the solver converges to the solution of the truncated system, and truncation of the operator and data sets the distance between that solution and the exact one. Although $\varepsilon$ does not bound the global error, it thus acts as an effective accuracy parameter.

\subsection{High-Resolution Three-Dimensional Elasticity}
\label{sec:results-exp2}

\begin{figure}
  \centering
  \includegraphics[width=\linewidth]{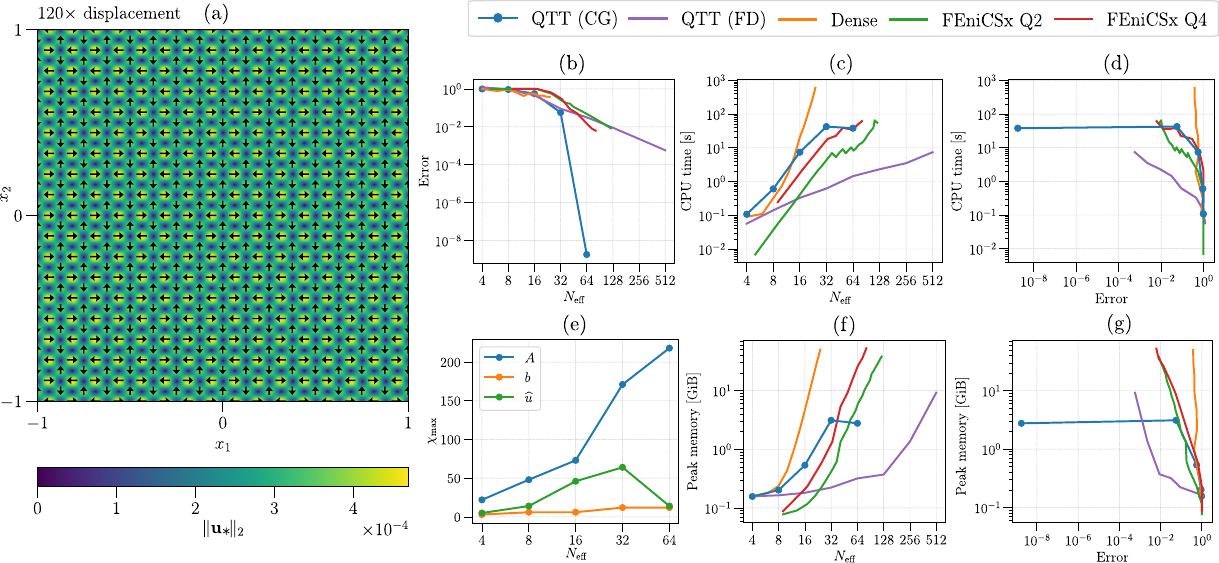}
  \caption{Three-dimensional benchmark with oscillatory load, comparing dense and QTT Chebyshev--Galerkin, QTT finite differences, and Q2 and Q4 finite elements, at the effective resolution $N_{\mathrm{eff}}$ defined in the text. (a) Displacement on the slice $x_3=0$; arrows show the displacement at local maxima of $|\bs{u}|$, magnified by a factor of 120. (b) Relative $\mathrm{L}^2$ error, (c) CPU time, and (f) peak memory versus $N_{\mathrm{eff}}$. (d) CPU time and (g) peak memory versus relative $\mathrm{L}^2$ error. (e) Maximum bond dimensions of the QTT operator $A$, load $b$, and solution $\hat{u}$.}
  \label{fig:exp2}
\end{figure}

We next compare the method with standard discretizations on a problem that requires high spatial resolution. The load is oscillatory, so that it must be finely resolved, but separable, so that it remains compressible. This isolates the effect of the discretization order from that of compression.

The domain is $\Omega=[-1,1]^3$, with a homogeneous material of $E=1$ and $\nu_{\mathrm P}=0.3$. The body force is a superposition of three separable oscillatory modes,
\begin{equation}
    \bs{f}(\x)=\sum_{i=1}^3 c_i\, \bs{\omega}_i(\x),
    \qquad
    \bs{\omega}_i(\x)=\sin(q_i x_i)\prod_{j\neq i}\cos(q_j x_j)\, \mathbf{e}_i,
    \label{eq:exp2-load}
\end{equation}
with wavenumbers $\bs{q}=\pi(8,9,10)$ and weights $\bs{c}=(1,-0.6,0.3)$, so that the shortest wavelength fits ten times across the domain [see Fig.~\ref{fig:exp2}(a)]. A bulk reaction $R_{\mu\nu}=\kappa\,\delta_{\mu\nu}$, with $\kappa=q_3^2$, keeps the reaction comparable to the elastic term. On each face $x_\gamma=\pm1$, the normal displacement $u_\gamma$ vanishes and the tangential tractions are zero, so that each face carries a Dirichlet condition on one component and a Neumann condition on the others. The exact solution is then $\bs{u}_*=\sum_i c_{*i}\,\bs{\omega}_i$, where the amplitudes $\bs{c}_*$ solve the $3\times3$ system
\begin{equation}
    \left[(\lambda_2|\bs{q}|^2+\kappa)\,I + (\lambda_1+\lambda_2)\,\bs{q}\bs{q}^{\intercal}\right]\bs{c}_* = \bs{c}.
\end{equation}

We compare dense and QTT Chebyshev--Galerkin discretizations with $N\leq24$ and $N\leq64$ modes per direction, respectively, with a QTT finite-difference (QTT-FD) discretization on uniform grids of up to $N=512$ points per direction, and with Q2 and Q4 finite elements in FEniCSx~\cite{alnaes2015,baratta2023}. To compare the methods at equal number of unknowns, we use the effective resolution $N_{\mathrm{eff}}$, equal to $N$ for the spectral and finite-difference methods and to $1+pN_{\mathrm{el}}$ for elements of degree $p$ with $N_{\mathrm{el}}$ elements per direction. Timings use all available threads, and FEniCSx uses the fastest MPI configuration at each resolution.

The results in Fig.~\ref{fig:exp2} show that spectral convergence ultimately outweighs the higher cost per unknown of the spectral operator. Once the oscillations are resolved, the spectral error drops rapidly and falls below $10^{-8}$ at $N=64$, whereas QTT-FD and finite elements converge algebraically and remain above $10^{-4}$ over the tested range [see Fig.~\ref{fig:exp2}(b)]. QTT-FD has the simplest operator and is the cheapest method at moderate accuracy, but its cost grows quickly with the requested accuracy. The error-based views in Figs.~\ref{fig:exp2}(d) and~(g) therefore show a crossover, beyond which spectral QTT reaches errors that the other methods do not attain within the tested time and memory.

The bond dimensions in Fig.~\ref{fig:exp2}(e) explain these costs. The load remains highly compressed, with $\chi_b\leq12$, whereas the bond dimension of the operator grows with $N$ and dominates the cost of spectral QTT. The bond dimension of the solution first increases and then drops sharply at $N=64$. Before the oscillations are resolved, the computed coefficients are an under-resolved approximation without a clear structure, and they compress poorly. Once the solution is resolved, its coefficients converge to those of the exact separable field, and the QTT recovers this structure.

\FloatBarrier

\subsection{High-Dimensional Parametric Cantilever Design}
\label{sec:results-exp3}

\begin{figure}
    \centering
    \includegraphics[width=0.7\linewidth]{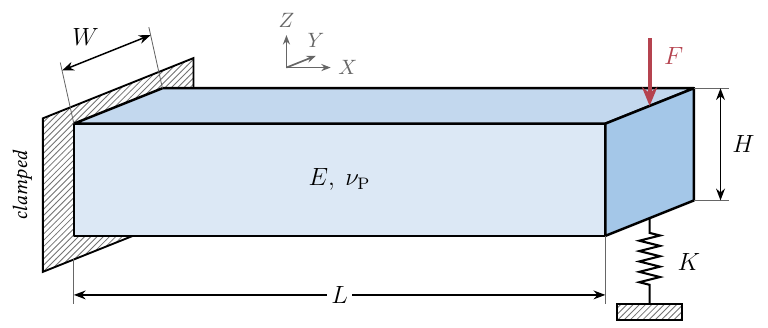}
    \caption{Parametric cantilever. The beam is clamped at $X=0$ and loaded at $X=L$ by a fixed transverse force $F$ and an elastic support of stiffness $K$. The length $L$, height $H$, width $W$, elastic modulus $E$, Poisson ratio $\nu_{\mathrm P}$, and stiffness $K$ are parameters.}
    \label{fig:exp3-model}
\end{figure}

\begin{figure}
    \centering
    \includegraphics[width=\linewidth]{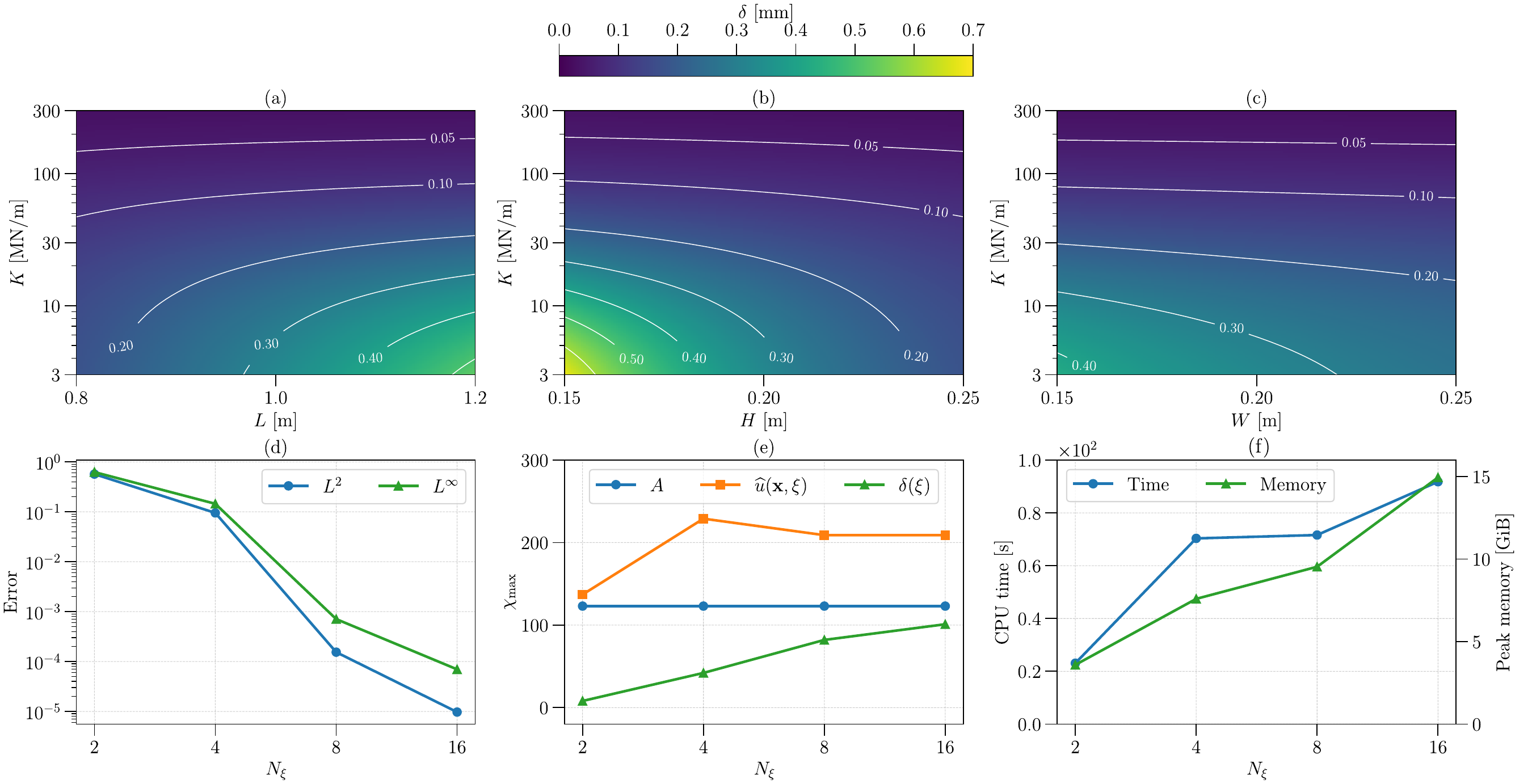}
    \caption{Parametric cantilever results. (a)--(c) Tip response $\delta(\bs{\xi})$ versus support stiffness $K$ and (a) length $L$, (b) height $H$, and (c) width $W$, with the remaining parameters at their nominal values, evaluated from the response QTT at $N_\xi=16$. (d) Relative $\mathrm{L}^2$ and $\mathrm{L}^\infty$ errors of the response over 128 random designs, versus $N_\xi$. (e) Maximum bond dimensions of the operator $A$, the displacement $\hat{u}(\x,\bs{\xi})$, and the response $\delta(\bs{\xi})$. (f) CPU time and peak memory.}
    \label{fig:exp3}
\end{figure}

We now turn to a design problem, in which the response of a structure is needed over a whole range of geometries, materials, and supports. We show that a single solve returns the response of every design, that its cost remains moderate as the number of designs grows by orders of magnitude, and that its error is small and decreases rapidly with the parametric resolution.

The structure is a three-dimensional cantilever of length $L$, width $W$, and height $H$ along $X$, $Y$, and $Z$, respectively (see Fig.~\ref{fig:exp3-model}). It is clamped at $X=0$, and at $X=L$ it carries a transverse load $F=10$~kN in the $-Z$ direction and rests on an elastic support of stiffness $K$, both distributed over the tip face,
\begin{equation}
    S_{\mu\nu} = \frac{K}{HW}\, \delta_{\mu3} \delta_{\nu3},
    \qquad
    g_\mu = -\frac{9F}{4HW}\left(1-x_2^2\right)\left(1-x_3^2\right)\delta_{\mu3},
    \label{eq:exp3-robin}
\end{equation}
so that they integrate to $K$ and $F$. The remaining faces are traction free. Six parameters $\bs{\xi}\in[-1,1]^6$ set
\begin{equation}
    \begin{gathered}
    L = 1 + 0.2\,\xi_1~\mathrm{m},
    \qquad
    H = 0.2 + 0.05\,\xi_2~\mathrm{m},
    \qquad
    W = 0.2 + 0.05\,\xi_3~\mathrm{m},
    \\
    E = 70 + 10\,\xi_4~\mathrm{GPa},
    \qquad
    \nu_{\mathrm P} = 0.3 + 0.05\,\xi_5,
    \qquad
    K = 30\times10^{\xi_6}~\mathrm{MN/m},
    \end{gathered}
    \label{eq:exp3-parameters}
\end{equation}
so that $K$ ranges from a weak to an almost rigid support, between 3 and 300~MN/m. The beam is mapped from $\Omega=[-1,1]^3$ by scaling each coordinate with $L$, $W$, and $H$ (Appendix~\ref{app:mappings}). The geometric parameters therefore enter the coefficients of the problem, and the operator depends on all six parameters. The quantity of interest is the load-weighted tip displacement
\begin{equation}
    \delta(\bs{\xi})=\frac{1}{F} \int_{X=L} g_\mu\, u_\mu(\x;\bs{\xi})\,\dS,
    \label{eq:exp3-response}
\end{equation}
which is obtained from the solution by contracting its spatial cores with the load, and is itself a QTT over the parameters.

We use $N=(16,8,8)$ spatial modes and $N_\xi=2,4,8,16$ modes per parameter, with $\varepsilon=10^{-10}$ and a residual tolerance of $10^{-5}$. To measure the error, we draw 128 designs at random among the $N_\xi^6$ Chebyshev--Gauss nodes of the parameter space and solve each of them independently with the dense discretization at the same spatial resolution. The error thus measures the parametric approximation, truncation, and solver tolerance, but not the spatial discretization.

The results in Fig.~\ref{fig:exp3} support the three claims. First, a single solve covers the whole design space. At $N_\xi=16$, the solution spans $16^6\approx1.7\times10^7$ designs and contains $5.15\times10^{10}$ coefficients, which would occupy 384~GiB as an explicit array, but its QTT occupies only 6.07~MiB, a compression factor of $6.5\times10^4$. The response $\delta(\bs{\xi})$ is a Chebyshev interpolant over the parameters, with bond dimension at most 101, and can be evaluated at any design at negligible cost. The maps of Figs.~\ref{fig:exp3}(a)--(c), with $161\times161$ designs each, are obtained in this way. They show, for instance, that the tip displacement is most sensitive to the height of the beam and that a stiff support suppresses the effect of its length.

Second, the cost grows only moderately with the number of designs. From $N_\xi=2$ to 16, the number of designs grows from 64 to $1.7\times10^7$, while the solve time grows from 23 to 92~s and the peak memory from 3.6 to 15~GiB [see Fig.~\ref{fig:exp3}(f)]. The bond dimension of the operator stays at 123, since refining the parameters does not introduce new couplings, and that of the solution stays below 230 [see Fig.~\ref{fig:exp3}(e)]. For comparison, a dense solve of a single design takes about 1.5~s on the same hardware, so a sweep over the same $1.7\times10^7$ designs would take about $2.5\times10^7$~s, roughly $3\times10^5$ times longer.

Third, the error is small and decreases rapidly with the parametric resolution [see Fig.~\ref{fig:exp3}(d)]. The relative $\mathrm{L}^2$ error over the sampled designs falls from $0.57$ at $N_\xi=2$ to $1.0\times10^{-5}$ at $N_\xi=16$, and the $\mathrm{L}^\infty$ error falls to $7\times10^{-5}$. The response of the whole design space is thus obtained in one solve, with an accuracy comparable to that of the solver tolerance.

\FloatBarrier

\subsection{High-Dimensional Intrusive Stochastic Elasticity}
\label{sec:results-exp4}

\begin{figure}
    \centering
    \includegraphics[width=0.45\linewidth]{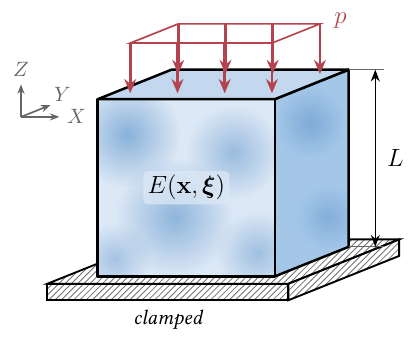}
    \caption{Stochastic block. The cube is clamped at its base and compressed by a uniform pressure $p$ on its top face. Its elastic modulus $E(\x,\bs{\xi})$ is a random field.}
    \label{fig:exp4-model}
\end{figure}

We finally consider uncertainty quantification, in which a material property is not known exactly but described statistically, and the goal is to compute the statistics of the response. We show that a single solve returns the complete random displacement field, from which means and variances follow directly; that its cost grows moderately with the number of random variables; and that it outperforms Monte Carlo sampling by orders of magnitude at equal accuracy.

The structure is an aluminum cube of side $L=10$~cm, clamped at its base and compressed by a uniform pressure $p=10$~MPa on its top face, with traction-free lateral faces (see Fig.~\ref{fig:exp4-model}). Its elastic modulus is a random field with mean $E_0=70$~GPa, a standard deviation of 5\%, and a correlation length of 7.5~cm, while $\nu_{\mathrm P}=0.3$. As is standard, we represent this field by a truncated Karhunen--Loève expansion~\cite{ghanem2003},
\begin{equation}
    E(\x; \bs{\xi})
    = E_0 \left[ 1 + \sum_{r=1}^M a_r\, \phi_r(\x)\, \xi_r\right],
    \label{eq:E_KL}
\end{equation}
in which smooth spatial modes $\phi_r$ are weighted by independent random variables $\xi_r$, uniformly distributed on $[-1,1]$. The construction of the modes and amplitudes is given in Appendix~\ref{app:kl}. With $M=10$ modes, the expansion captures 99.3\% of the variance of the material and keeps $E\geq0.75E_0$ for every realization. The random variables enter as parameters, as in Sec.~\ref{sec:qtt-parametric}. For uniform variables, integrating the weak form over $\Xi$ is proportional to taking its expectation, so this is the intrusive stochastic Galerkin method~\cite{ghanem2003,matthies2005}.

The quantity of interest is the average compression of the top face,
\begin{equation}
    \delta(\bs{\xi}) = -\frac{1}{L^2} \int_{Z=L} u_3(\x;\bs{\xi})\,\dS,
    \label{eq:exp4-response}
\end{equation}
which, as for the cantilever, is obtained by contracting the spatial cores of the solution and is itself a QTT over the random variables.

We vary the number of random variables from $M=1$ to 10, with $N_\xi=2,4,8$ modes per variable and $N=8$ spatial modes per direction, with $\varepsilon=10^{-10}$ and a residual tolerance of $10^{-5}$. At $M=10$ and $N_\xi=8$, the solution contains $1.65\times10^{12}$ coefficients. To measure the error, we draw 100 realizations of $\bs{\xi}$ uniformly at random, evaluate the QTT solution at them by Chebyshev interpolation, and compare it with independent dense solves at the same spatial resolution.

\begin{figure}
    \centering
    \includegraphics[width=\linewidth]{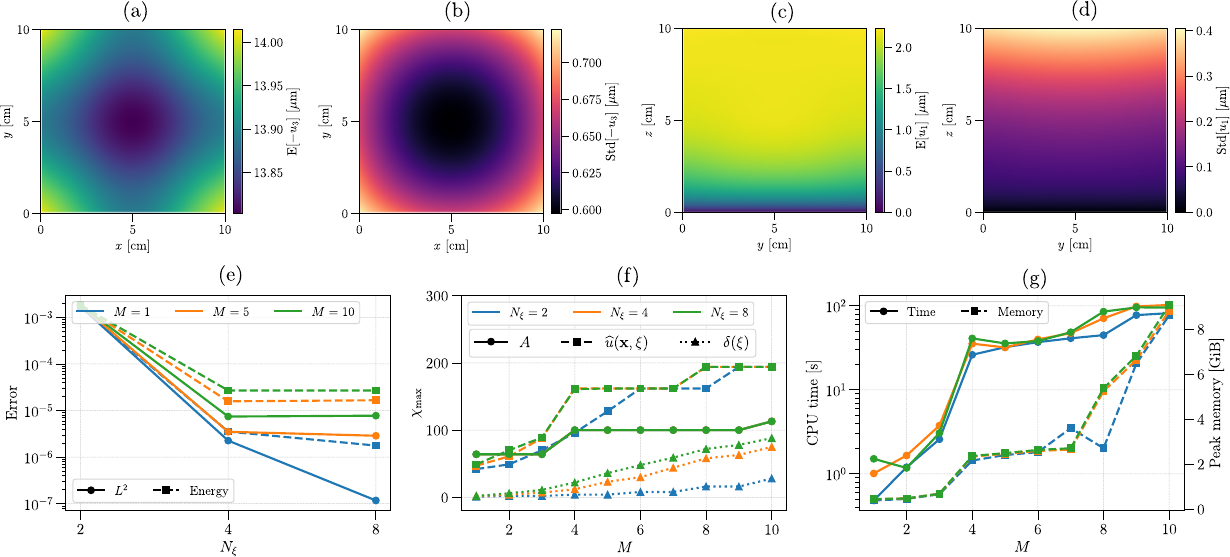}
    \caption{Stochastic block results. (a, b) Mean and standard deviation of the vertical displacement on the top face, and (c, d) of the outward displacement on the lateral face $X=L$, for $M=10$ and $N_\xi=8$, obtained from the QTT without sampling error. (e) Relative $\mathrm{L}^2$ and energy-norm errors over 100 random realizations versus $N_\xi$, for $M=1,5,10$. (f) Maximum bond dimensions of the operator $A$, the displacement $\hat{u}(\x,\bs{\xi})$, and the response $\delta(\bs{\xi})$ versus $M$. (g) CPU time and peak memory versus $M$.}
    \label{fig:exp4}
\end{figure}

The results in Fig.~\ref{fig:exp4} support the three claims. First, a single solve returns the full random field, from which any statistic follows by integrating over the parameter cores, without sampling. The mean vertical displacement of the top face is about $14~\mu$m, with a standard deviation between 0.60 and $0.72~\mu$m that is smallest near the center, and the lateral faces bulge outward by up to $2.2~\mu$m [see Figs.~\ref{fig:exp4}(a)--(d)]. The solution is accurate. From $N_\xi=4$ on, the relative $\mathrm{L}^2$ error over the sampled realizations stays below $2\times10^{-5}$ for every $M$, and at $M=10$ and $N_\xi=8$ it is $7.8\times10^{-6}$, with an energy-norm error of $2.7\times10^{-5}$, consistent with the residual tolerance [see Fig.~\ref{fig:exp4}(e)].

Second, the cost grows moderately with the number of random variables. From $M=1$ to 10, the number of coefficients grows by a factor of $10^8$, while the solve time grows from 1.5 to 95~s and the peak memory from 0.45 to 9.1~GiB [see Fig.~\ref{fig:exp4}(g)]. The bond dimension of the operator only grows from 64 to 113, and that of the solution from 49 to 194 [see Fig.~\ref{fig:exp4}(f)]. At $M=10$, the QTT of the solution occupies 5.18~MiB, against 12~TiB for the explicit array, a compression factor of $2.4\times10^6$.

Third, the method outperforms Monte Carlo sampling. A dense solve of one realization takes about 0.025~s on the same hardware, so the time of the QTT solve allows about 3800 Monte Carlo samples. These estimate the mean and standard deviation of $\delta$ in Eq.~\eqref{eq:exp4-response} with relative errors of about $7\times10^{-4}$ and $10^{-2}$, respectively, against errors of order $10^{-5}$ for the QTT solution. Reaching an error of $10^{-5}$ in the mean would require about $1.7\times10^{7}$ samples, roughly 5 days of computation. This comparison is specific to the present observable, for which a single random variable already captures 99.5\% of the response variance (see Appendix~\ref{app:kl}), so that methods exploiting this structure, such as quasi-Monte Carlo~\cite{caflisch1998} or sparse polynomial chaos and stochastic collocation~\cite{blatman2011,babuvska2007}, may narrow the gap. The case $M=10$ serves instead to test the method in a high-dimensional setting.

\section{Discussion}
\label{sec:discussion}

We have introduced a QTT Chebyshev--Galerkin method that assembles and solves global spectral systems directly in compressed form. Starting from the weak formulation, every ingredient of the discrete problem, including variable coefficients, coupled fields, mapped geometries, Robin terms, and Dirichlet constraints, is built from a few one-dimensional Chebyshev blocks through Kronecker products, sums, and MPO products. Boundary conditions are thus imposed systematically and in a structure-preserving way, by operations on whole operators rather than by editing individual entries, and the resulting SPD systems are solved by variational tensor-network algorithms. Parameters and random inputs enter as additional coordinates, so the same construction returns the solution over an entire parameter space in a single solve.

The experiments support two conclusions. First, compression does not compromise spectral accuracy. In the manufactured problem, the QTT solution follows the convergence of the dense discretization until it reaches the truncation floor, with solution bond dimensions of at most 19. In the three-dimensional benchmark, it reaches errors below $10^{-8}$ at $N=64$, whereas QTT finite differences and Q2/Q4 finite elements remain above $10^{-4}$ over the tested range, although QTT finite differences are cheaper at moderate accuracy. Second, the cost is governed by the bond dimensions of the operator and the solution rather than by the number of unknowns. A single solve of under 100~s covers $1.7\times10^7$ cantilever designs, with $5.2\times10^{10}$ nominal coefficients stored in 6~MiB, and a stochastic problem with ten random coordinates and $1.65\times10^{12}$ nominal coefficients is stored in 5~MiB, a compression factor of $2.4\times10^6$. In both cases, the error against independent dense solves is around $10^{-5}$.

The experiments also show the limits of the approach. Spectral accuracy requires smooth solutions and smooth mappings, and compression requires low-rank structure. Moreover, the cost depends on the operator as much as on the solution: in the three-dimensional benchmark, the bond dimension of the operator grows with resolution even when the solution remains highly compressible. The ill-conditioning of the spectral stiffness operator further limits the attainable accuracy, since truncation errors are amplified by the condition number, and is the main practical bottleneck of the current implementation.

These limitations point to the natural extensions. Better conditioning could be obtained from sparse, well-conditioned spectral bases~\cite{olver2013} or from multilevel preconditioners built directly in QTT form. Nonlinear problems, such as large-deformation elasticity, can be addressed by Newton or fixed-point iterations, in which each step is a linear problem of the form treated here, with the current solution entering through the multiplication operator. Time-dependent problems can be treated either by implicit time stepping, in which each step is again an elliptic Galerkin problem, or by treating time as an additional Chebyshev coordinate in a space--time formulation~\cite{dolgov2012parabolic}. Finally, the method is best suited to problems with smooth structure in many coordinates, where explicit discretizations or sampling become prohibitive. Intrusive uncertainty quantification and mechanical design, as illustrated here, are natural candidates, together with optimization and inverse problems that require repeated access to a parametric solution. 

\section*{Acknowledgements}

We thank Jorge Gidi for implementing the QTT finite-difference method used as a benchmark in this work, and for many insightful discussions. This work received funding from the EuroHPC Joint Undertaking under the QEC4QEA project (Grant No. 101194322). JJGR and JJRA acknowledge support from the CSIC Interdisciplinary Thematic Platform on Quantum Technologies (PTI-QTEP+). This work was also supported by the Spanish Ministry of Science, Innovation and Universities and the State Research Agency (MICIU/AEI/10.13039/501100011033) through Project No. PDC2022-133486-I00.







\appendix

\section{Smooth Coordinate Mappings}
\label{app:mappings}

The reference domain $\Omega = [-1, 1]^d$ can be mapped to a curved physical domain $\wt{\Omega}$ through a smooth, orientation-preserving diffeomorphism $\Phi:\Omega\to\wt{\Omega}$, with $\wt{\x} = \Phi(\x)$. Its Jacobian matrix and determinant are
\begin{equation}
    F_{\alpha' \alpha}(\x)=\frac{\partial \wt{x}_{\alpha'}}{\partial x_\alpha},
    \qquad
    J(\x)=\det F(\x)>0.
\end{equation}
Changing variables in the weak form of Eq.~\eqref{eq:bvp_weak} yields a problem of the same form on $\Omega$, with coefficients
\begin{align}
    \mathbb{C}_{\mu\alpha\nu\beta}(\x)
    &= J(\x)\, F^{-1}_{\alpha \alpha'}(\x) \,
    \wt{\mathbb{C}}_{\mu\alpha'\nu\beta'}\!\left(\Phi(\x)\right) F^{-1}_{\beta \beta'}(\x),
    \\[0.5em]
    R_{\mu\nu}(\x) &=J(\x) \, \wt{R}_{\mu\nu}(\Phi(\x)),
    \qquad
    f_\mu(\x) = J(\x) \, \wt{f}_\mu(\Phi(\x)).
\end{align}
On a face of $\Omega$ with outward normal $\bs{n}$, the boundary data are scaled by the surface Jacobian $J_\Gamma(\x)=J(\x)\,|F^{-\intercal}(\x)\,\bs{n}|$,
\begin{equation}
    S_{\mu\nu}(\x) = J_\Gamma(\x)\,\wt{S}_{\mu\nu}(\Phi(\x)),
    \qquad
    g_\mu(\x) = J_\Gamma(\x)\,\wt{g}_\mu(\Phi(\x)),
\end{equation}
while Dirichlet data are simply composed with the mapping, $\bar{u}_\mu(\x)=\wt{\bar{u}}_\mu(\Phi(\x))$. Only the derivative indices $\alpha,\beta$ are transformed. The field components $u_\mu$ remain the Cartesian components of the physical field, so for elasticity the solution is still the physical displacement. Since $J>0$ and $F$ is invertible, the pulled-back problem inherits the symmetry and coercivity of the original one.

The geometry is thus absorbed into spatially varying coefficients on $\Omega$, which are loaded and multiplied like any other data (Sec.~\ref{sec:qtt-chebyshev}). Its only cost is the bond dimension of these coefficients, which remains small for smooth mappings. The mapping may also depend on the parameters, $\Phi=\Phi(\x,\bs{\xi})$, as in the parametric cantilever of Sec.~\ref{sec:results-exp3}. Throughout this work, we assume that any such mapping has already been incorporated into the problem data.

\section{QTT Chebyshev Transform}
\label{app:transform}

We construct the one-dimensional transform $\mathcal{D}_\gamma$ of Eq.~\eqref{eq:chebyshev_transform} for $N=2^n$ nodes as a product of $2n-1$ exact MPO layers, each with bond dimension at most 10, independently of $N$.

Up to a diagonal rescaling, $\mathcal{D}_\gamma$ is the orthonormal DCT-II,
\begin{equation}
    \mathcal{D}_\gamma = S_N\, C_N^{\rm{II}},
    \qquad
    \Big[C_N^{\rm{II}}\Big]_{kj}
    =\sqrt{\frac{2-\delta_{k0}}{N}} \cos\!\left(\frac{k(2j+1)\pi}{2N}\right),
    \qquad
    \Big[S_N\Big]_{kk} = \sqrt{\frac{2 - \delta_{k0}}{N}},
\end{equation}
for $k,j=0,\ldots,N-1$. The rescaling $S_N$ is an exact MPO of bond dimension 2, so it remains to represent $C_N^{\rm{II}}$. Like the DFT, the DCT-II matrix admits a recursion that halves the transform size at each step~\cite{perera2016},
\begin{equation}
    C_M^{\rm{II}} = P_M^\intercal
    \begin{pmatrix}
        C_{M/2}^{\rm{II}} & 0 \\ 0 & C_{M/2}^{\rm{IV}}
    \end{pmatrix}
    B_M,
    \qquad
    C_M^{\rm{IV}} = P_M^\intercal U_M
    \begin{pmatrix}
        C_{M/2}^{\rm{II}} & 0 \\ 0 & C_{M/2}^{\rm{II}}
    \end{pmatrix}
    R_M,
\end{equation}
where $C_M^{\rm{IV}}$ is the orthonormal DCT-IV matrix, $[C_M^{\rm{IV}}]_{kj}=\sqrt{2/M}\cos[(2k+1)(2j+1)\pi/(4M)]$, $P_M$ is the even--odd permutation, $B_M$ is a butterfly matrix, and $R_M$, $U_M$ are sparse rotation factors. Since $M=2^m$, each halving step acts on a single binary digit of the index. Collecting the operations at the same recursion depth therefore gives one MPO layer per digit. Which branch, DCT-II or DCT-IV, applies at a given depth depends only on the digits already processed, and this information is carried by the bond indices, which bounds every layer at bond dimension 10. This yields $n$ butterfly layers and $n-2$ correction layers from the DCT-IV branches, which together with $S_N$ give $2n-1$ layers in total.

\begin{figure}
    \centering
    \includegraphics[width=0.8\linewidth]{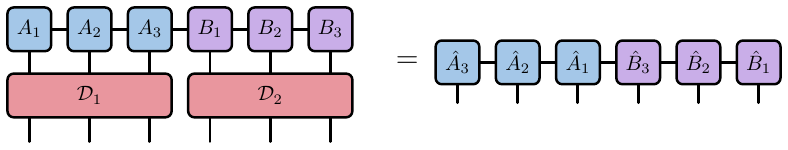}
    \caption{Multidimensional QTT Chebyshev transform applied independently to two tensor factors,
    $A\otimes B \xrightarrow{\mathcal D} \wt{A}\otimes\wt{B}$, each represented here by $n=3$ binary QTT cores. The radix-two transform acts locally on each coordinate block and directly produces bit-reversed modal indices, thus defining the native spectral QTT layout.}
    \label{fig:transform}
\end{figure}

As in the FFT, the recursion produces the output digits in reversed order of significance. Restoring the natural order would require a bit-reversal permutation, which is highly nonlocal in QTT form~\cite{chen2023}. We therefore keep the Chebyshev coefficients in this bit-reversed layout and build the spectral blocks of Appendix~\ref{app:mpos} directly in it, which only requires reversing the order of their cores. For multivariate data, the transform acts independently on the cores of each coordinate (see Fig.~\ref{fig:transform}).

\section{MPO Constructions for the Chebyshev Blocks}
\label{app:mpos}

This appendix gives the MPOs of the three one-dimensional blocks of Sec.~\ref{sec:qtt-chebyshev}. Throughout, $N=2^n$, and the binary digits of the row and column indices $i$ and $j$ are assigned to the MPO cores from most to least significant. The bit-reversed layout of Appendix~\ref{app:transform} is then obtained by reversing the order of the cores.

\subsection{Differentiation}
\label{app:mpos-differentiation}

The differentiation matrix of Sec.~\ref{sec:qtt-chebyshev} can be written compactly as
\begin{equation}
  D^{(1)}_{ij} = \frac{2j}{1+\delta_{i0}}\, \mathbf{1}_{j>i}\, \mathbf{1}_{i+j\ {\rm odd}}.
  \label{eq:chebyshev_derivative_entries}
\end{equation}
Its MPO reads the digit pairs $(i_p,j_p)$, $p=1,\ldots,n$, one at a time. Six bond states suffice to record whether $j>i$ has already been decided and whether $i=0$, while the factor $2j$ is accumulated digit by digit. For fixed $i_p,j_p\in\{0,1\}$, the core is
\begin{equation}
    D^{[p]}_{i_p j_p}
=
\begin{pmatrix}
    \delta_{i_pj_p} & 0 & j_p(1-i_p) & w_p j_p(1-i_p) & 0 & w_p i_pj_p
    \\
    0 & 1-i_p & 0 & 0 & -\tfrac{1}{2}w_p j_p(1-i_p) & 0
    \\
    0 & 0 & 1 & w_p j_p & 0 & 0
    \\
    0 & 0 & 0 & 1 & 0 & 0
    \\
    0 & 0 & 0 & 0 & 1-i_p & 0
    \\
    0 & 0 & 0 & j_p(1-i_p) & 0 & \delta_{i_pj_p}
\end{pmatrix},
\qquad w_p = 2^{n-p+1},
\end{equation}
where $w_p$ is the contribution of $j_p=1$ to $2j$. The chain is closed by the boundary vectors $(1,1,0,0,0,0)$ on the left and $(0,0,0,1,1,0)^\intercal$ on the right. The parity condition is imposed on the last, least-significant digit pair by multiplying the last core by $\mathbf{1}_{i_n\neq j_n}=\mathbf{1}_{i+j\ \rm odd}$. The result is exact, with bond dimension 6 for every $N$.

\subsection{Gram matrix}
\label{app:mpos-gram}

The Gram and multiplication matrices are both built from the same four structured operators. Given a vector $v$, with $v_k=0$ outside its range, we define
\begin{equation}
    \mathsf{H}(v)_{ij} = v_{i+j},
    \qquad
    \mathsf{L}(v)_{ij} = v_{i-j}\,\mathbf{1}_{i\geq j},
    \qquad
    \mathsf{U}(v)_{ij} = v_{j-i}\,\mathbf{1}_{i\leq j},
    \qquad
    \mathsf{S}(v)_{ij} = \delta_{i0}\,v_j,
    \label{eq:structured_operators}
\end{equation}
namely a Hankel matrix, the lower and upper triangular parts of a Toeplitz matrix, and a matrix whose only nonzero row is the first. With $w_k=2/(1-k^2)$ for even $k$ and $w_k=0$ for odd $k$, Eq.~\eqref{eq:gram_chebyshev} becomes $G^{(1)}_{ij}=\tfrac12(w_{i+j}+w_{|i-j|})$, that is,
\begin{equation}
    G^{(1)} = \tfrac12\big[\mathsf{H}(w)+\mathsf{L}(w)+\mathsf{U}(w)\big] - I,
\end{equation}
where the identity removes the diagonal counted by both $\mathsf{L}$ and $\mathsf{U}$, since $w_0=2$.

Each operator in Eq.~\eqref{eq:structured_operators} is obtained directly from the QTT of $v$, core by core. Every core of $v$ becomes an MPO core acting on the digit pair $(i_p,j_p)$, and its bond index is extended by one binary flag. For $\mathsf{H}$, the flag is the carry of the binary sum $i+j$. For $\mathsf{L}$ and $\mathsf{U}$, it is the borrow of the difference $i-j$, which also decides the sign condition. Entries whose index falls outside the range of $v$ are discarded at the end of the chain. The flag doubles the bond dimension, so $\mathsf{H}$, $\mathsf{L}$, and $\mathsf{U}$ have bond dimension $2\chi_v$, while $\mathsf{S}$ only fixes $i=0$ and keeps $\chi_v$. For the Gram matrix, this gives $\chi_G\leq 6\chi_w+1$.

The construction involves no matrix products and no quadrature. Its only approximation is the TT-SVD of the vector $w$ of length $2N$, whose bond dimension is $\chi_w=9$ at machine precision for all tested resolutions up to $N=2^{20}$. The resulting MPO therefore has bond dimension at most $\chi_G=55$, independently of $N$. Since this bound is not tight, a final compression at machine precision reduces it further to $\chi_G\approx 23$, which saturates from $N=2^{8}$ onwards.

\subsection{Multiplication matrix}
\label{app:mpos-multiplication}

In terms of the operators of Eq.~\eqref{eq:structured_operators}, Eq.~\eqref{eq:Q_chebyshev} reads
\begin{equation}
    Q(a) = \tfrac12\big[\mathsf{H}(\hat a)+\mathsf{L}(\hat a)+\mathsf{U}(\hat a)-\mathsf{S}(\hat a)\big].
\end{equation}
Here, $\mathsf{L}$ accounts for the term $\delta_{i\,j+k}$, $\mathsf{H}$ for $\delta_{i\,|k-j|}$ with $k\geq j$, and $\mathsf{U}$ for $\delta_{i\,|k-j|}$ with $k\leq j$. The term $\mathsf{S}$ removes the first row, where $k=j$ is counted by both $\mathsf{H}$ and $\mathsf{U}$. The same construction as for the Gram matrix therefore gives $Q(a)$ exactly from the QTT of $\hat a$, with $\chi_Q\leq7\chi_{\hat a}$. For multivariate coefficients, it is applied to the cores of each coordinate in turn, so $Q(a)$ is obtained directly from the QTT of $\hat{a}_{\bs{k}}$ without any separability assumption.

\section{Preconditioning of the QTT Linear Systems}
\label{app:preconditioning}

The Chebyshev differentiation matrix is ill-conditioned: its entries grow with the spectral degree, so the eigenvalues of the differential operator span several orders of magnitude. Moreover, since the Chebyshev basis is not orthogonal under the unweighted Galerkin measure, the Gram matrices are themselves ill-conditioned, with $\kappa(G^{(1)})=\O{N}$. Since the condition number of a Kronecker product is the product of those of its factors, this effect compounds multiplicatively, $\kappa(A\otimes B)=\kappa(A)\,\kappa(B)$. For instance, following Eqs.~\eqref{eq:A_bulk} and~\eqref{eq:galerkin_system_lifted}, the operator of $-\Delta u+u$ on $[-1,1]^d$ with homogeneous Dirichlet conditions reads
\begin{equation}
    A = Z^\intercal\Big(\sum_{\alpha=1}^d D_\alpha^\intercal G\,D_\alpha + G\Big)Z,
    \qquad
    D_\alpha^\intercal G\,D_\alpha = G^{(1)}_1\otimes\cdots\otimes D^{(1)\intercal}_\alpha G^{(1)}_\alpha D^{(1)}_\alpha\otimes\cdots\otimes G^{(1)}_d,
    \label{eq:model_operator}
\end{equation}
where the second identity follows from Eq.~\eqref{eq:DG_kronecker}. Each derivative term thus carries the Gram matrices of all the other directions, and each parametric coordinate contributes a further Gram factor $G^{(1)}_{\xi_q}$. As a result, the condition number grows exponentially with the number of spatial and parametric coordinates (Table~\ref{tab:conditioning}, first column).

\begin{table}
    \centering
    \caption{Condition number of $-\nabla\cdot(a\nabla u)+a\,u$ on $[-1,1]^d$ with homogeneous Dirichlet conditions and $p$ parameters, for $N=16$, $N_\xi=8$, and $a=1+\frac{1}{2p}\sum_q\xi_q$ ($a=1$ for $p=0$). Jacobi denotes the global diagonal scaling $\operatorname{diag}(A)^{-1/2}A\operatorname{diag}(A)^{-1/2}$, and Cholesky the change of basis with $\ell_\gamma=1$ and $\tau=100$.}
    \label{tab:conditioning}
    \begin{tabular}{cccccc}
        \toprule
        $d$ & $p$ & None & Jacobi & Cholesky & Cholesky + Jacobi \\
        \midrule
        1 & 0 & $6.6\times10^{1}$ & $3.8$ & $2.8\times10^{1}$ & $2.2$ \\
        2 & 0 & $7.6\times10^{2}$ & $3.7\times10^{2}$ & $1.6\times10^{1}$ & $5.6$ \\
        3 & 0 & $6.3\times10^{4}$ & $3.6\times10^{4}$ & $1.3\times10^{2}$ & $4.4\times10^{1}$ \\
        \midrule
        1 & 1 & $1.3\times10^{3}$ & $5.2\times10^{1}$ & $8.0\times10^{1}$ & $6.2$ \\
        1 & 2 & $1.2\times10^{4}$ & $3.3\times10^{2}$ & $8.0\times10^{1}$ & $6.2$ \\
        1 & 3 & $1.3\times10^{5}$ & $2.1\times10^{3}$ & $8.0\times10^{1}$ & $6.2$ \\
        \bottomrule
    \end{tabular}
\end{table}

We address these two effects with two complementary techniques. Jacobi preconditioning acts locally: it rescales the diagonal of the systems solved by AMEn and thereby compensates the disparity of scales introduced by $D^{(1)}$. A Cholesky-based change of basis acts globally: it orthonormalizes the basis along each coordinate and thereby breaks the multiplicative growth. Together, they reduce the condition number by up to four orders of magnitude (Table~\ref{tab:conditioning}, last column). In practice, we adapt the preconditioning to the difficulty of the problem: for the purely spatial problems of Secs.~\ref{sec:results-exp1} and~\ref{sec:results-exp2}, Jacobi preconditioning alone suffices, whereas for the high-dimensional parametric problems of Secs.~\ref{sec:results-exp3} and~\ref{sec:results-exp4}, we combine it with the Cholesky change of basis.

\subsection{Jacobi preconditioning}
At each step, AMEn solves a small linear system $A_kx_k=b_k$ for the active core~\cite{dolgov2014}. We solve it by conjugate gradients preconditioned with $\operatorname{diag}(A_k)^{-1}$~\cite{rohrig2025}, whose entries are obtained by contracting the diagonal slices of the MPO core and its environments, without forming $A_k$. Rescaling by the diagonal brings all modes to a comparable scale; in one dimension, it reduces the growth of $\kappa$ from $\O{N^2}$ to $\O{N}$. However, it cannot correct the nonorthogonality of the basis, so on its own it barely helps in several coordinates (Table~\ref{tab:conditioning}, second column). For the benchmark of Sec.~\ref{sec:results-exp2}, we use a block-Jacobi variant that retains the couplings within groups of local unknowns.

\subsection{Cholesky change of basis}
To break the multiplicative growth, we apply a second change of basis on top of the recombination $Z$ of Sec.~\ref{sec:qtt-dirichlet}, one coordinate at a time. This amounts to a Kronecker-product preconditioner whose factors are shifted one-dimensional stiffness matrices, as introduced for the discrete Poisson problem in Ref.~\cite{vanloan1993}. Along each spatial coordinate, we take as reference the one-dimensional counterpart of Eq.~\eqref{eq:model_operator},
\begin{equation}
    H_\gamma = Z^{(1)\intercal}_\gamma\Big(\ell_\gamma^{-2}\,D^{(1)\intercal}_\gamma G^{(1)}_\gamma D^{(1)}_\gamma + \tau\,G^{(1)}_\gamma\Big)Z^{(1)}_\gamma,
\end{equation}
where $\ell_\gamma$ is the physical length of the domain and $\tau$ is a shift ($\tau=100$ in our experiments). Along a parametric coordinate, it reduces to the Gram matrix $H_\gamma=G^{(1)}_{\xi_q}$. We factorize each reference as $H_\gamma=U_\gamma^\intercal U_\gamma$ and use $\widetilde{Z}_\gamma=U_\gamma^{-1}$ as the new basis, in which $H_\gamma$ becomes the identity, $\widetilde{Z}_\gamma^\intercal H_\gamma \widetilde{Z}_\gamma=I$. Defining $\widetilde{Z}=\bigotimes_\gamma \widetilde{Z}_\gamma$, AMEn then solves the equivalent system
\begin{equation}
    \big(\widetilde{Z}^\intercal\Ahom \widetilde{Z}\big)\,y = \widetilde{Z}^\intercal\bhom,
    \qquad
    \hat{u} = \ulift + Z\widetilde{Z}\,y,
\end{equation}
so that $Z$ imposes the Dirichlet conditions and $\widetilde{Z}$ improves the conditioning. Along the parametric coordinates, the basis becomes orthonormal and the Gram factors in the operator become identities. The parameters then contribute only through the variation of the coefficients: if these depend on $\bs{\xi}$ through a factor bounded between $c_{\min}$ and $c_{\max}$, the operator is bounded by the same factors times its parameter-independent counterpart, so the condition number increases at most by a factor $c_{\max}/c_{\min}$, independently of $N_\xi$ and of the number of parameters (Table~\ref{tab:conditioning}, last three rows). Along the spatial coordinates, the situation is less favorable: since the operator in Eq.~\eqref{eq:model_operator} is a sum of Kronecker products, $\widetilde{Z}$ cannot normalize all its terms at once, but the remaining spread is largely a scaling of the modes, which the Jacobi preconditioner then removes.

This combination has worked well for the problems considered here, but leaves room for improvement. In several spatial dimensions, the preconditioned condition number still grows with $N$, and variable coefficients, mappings, and couplings between components enter only through the diagonal. Sparse, well-conditioned spectral bases~\cite{olver2013} or multilevel preconditioners in QTT form~\cite{bachmayr2020stability} could make it independent of the resolution. Since conditioning currently limits the attainable accuracy, we regard this as an interesting direction for future work.

\section{Spectral Transfinite Interpolation}
\label{app:tfi}

The lift $\ulift$ of Sec.~\ref{sec:qtt-dirichlet} only needs to match the Dirichlet data on the Dirichlet faces, since the correction $Z\uhom$ accounts for the interior. We construct it by transfinite interpolation (TFI)~\cite{gordon1973}, directly in Chebyshev coefficient space, assuming that the data of adjacent faces agree where they meet.

In one dimension, values $\bar{u}_\pm$ at $x=\pm1$ are extended to the interval as $\bar{u}_-B^-(x)+\bar{u}_+B^+(x)$, where each blending function equals one at its own endpoint and vanishes at the other. We use the cubic Hermite functions
\begin{equation}
    B^{\pm}(x) = \tfrac14(1\pm x)^2(2\mp x) = \tfrac12 T_0(x) \pm \tfrac{9}{16}T_1(x) \mp \tfrac{1}{16}T_3(x),
    \label{eq:tfi_hermite_functions}
\end{equation}
with coefficient vectors $b^\pm$, so that $t_\tau^\intercal b^\sigma=\delta_{\sigma\tau}$ in terms of the boundary vectors of Eq.~\eqref{eq:trace_operator}. In several dimensions, data on the face $x_\gamma=\sigma$ are extended by multiplying them by $B^\sigma(x_\gamma)$, which in coefficient space is the counterpart of the trace onto that face,
\begin{equation}
    \mathsf{E}_{\gamma\sigma} = I\otimes\cdots\otimes b^\sigma\otimes\cdots\otimes I,
    \qquad
    \mathsf{T}_{\gamma\sigma} = I\otimes\cdots\otimes t_\sigma^\intercal\otimes\cdots\otimes I.
\end{equation}
The extension thus reproduces the data on its own face, $\mathsf{T}_{\gamma\sigma}\mathsf{E}_{\gamma\sigma}=I$, and vanishes on the opposite one, $\mathsf{T}_{\gamma\,-\sigma}\mathsf{E}_{\gamma\sigma}=0$.

\begin{figure}
    \centering
    \includegraphics[width=\linewidth]{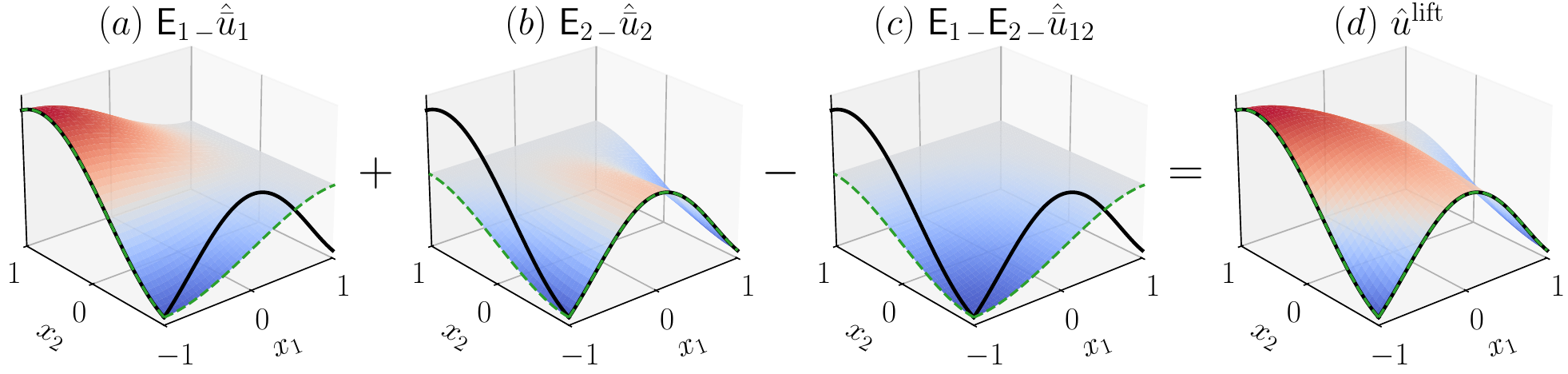}
    \caption{Transfinite interpolation on two adjacent Dirichlet faces, $x_1=-1$ and $x_2=-1$. Black lines show the data and dashed green lines the trace of each function on both faces. Each extension (a,b) matches the data only on its own face; subtracting the extended corner value (c) yields a lift (d) that matches both.}
    \label{fig:tfi}
\end{figure}

Summing the extensions of all faces is not enough, since they interfere where faces meet. For two adjacent faces $x_1=-1$ and $x_2=-1$ with data $\hat{\bar{u}}_1$ and $\hat{\bar{u}}_2$, the extension of each one adds the common corner value $\hat{\bar{u}}_{12}$ on the other face, and subtracting the extension of this value corrects both (see Fig.~\ref{fig:tfi}),
\begin{equation}
    \ulift = \mathsf{E}_{1\,-}\hat{\bar{u}}_1 + \mathsf{E}_{2\,-}\hat{\bar{u}}_2 - \mathsf{E}_{1\,-}\mathsf{E}_{2\,-}\hat{\bar{u}}_{12}.
\end{equation}
In general, let $\mathcal{S}_\mu$ be the set of Dirichlet faces of component $\mu$. Each subset $\mathcal{C}\subseteq\mathcal{S}_\mu$ with at most one face per coordinate meets on a face, edge, or corner of $\Omega$, whose data $\hat{\bar{u}}_{\mu,\mathcal{C}}$ follow by tracing the data of any face in $\mathcal{C}$ onto the others. The lift is then the inclusion--exclusion sum over all such nonempty subsets,
\begin{equation}
    \ulift_\mu
    =
    \sum_{\mathcal{C}}
    (-1)^{|\mathcal{C}|-1}
    \Bigg[
        \prod_{(\gamma,\sigma)\in \mathcal{C}}\mathsf{E}_{\gamma\sigma}
    \Bigg]
    \hat{\bar{u}}_{\mu,\mathcal{C}},
    \label{eq:spectral_tfi_lift}
\end{equation}
in which single faces extend their data, pairs subtract the edge data counted twice, triples add back the corner data, and so on. Each vector $b^\sigma$ has QTT bond dimension two, so each $\mathsf{E}_{\gamma\sigma}$ is an exact MPO of bond dimension two. The lift is thus assembled in QTT form from Kronecker products and sums alone, without evaluating it on any grid.

\section{Manufactured Solution of the Verification Problem}
\label{app:manufactured}

The manufactured displacement of Sec.~\ref{sec:results-exp1} is built from Euler--Bernoulli cantilever kinematics. Let
\begin{equation}
    \phi_{\mathrm{stat}}(s)=\frac{s^2(s^2-4s+6)}{3},
    \qquad
    \phi_3(s)=\frac{\cosh(\beta_3s)-\cos(\beta_3s)-\kappa_3\left[\sinh(\beta_3s)-\sin(\beta_3s)\right]}{\cosh\beta_3-\cos\beta_3-\kappa_3\left[\sinh\beta_3-\sin\beta_3\right]},
\end{equation}
be the normalized static deflection and the third clamped--free bending mode, with $\beta_3\approx7.85$ the third positive root of $\cos\beta\cosh\beta=-1$ and $\kappa_3=(\cosh\beta_3+\cos\beta_3)/(\sinh\beta_3+\sin\beta_3)$. With the deflection and centerline
\begin{equation}
    w(X)=-\delta\ell\left[\tfrac{7}{8}\phi_{\mathrm{stat}}(X/\ell)+\tfrac{1}{8}\phi_3(X/\ell)\right],
    \qquad
    Y_c(X)=\tfrac12\left[1+h(X/\ell)\right],
\end{equation}
and $\delta=2\times10^{-3}$, the displacement is manufactured as
\begin{equation}
    \bs{u}_*(X,Y)
    = \begin{pmatrix}
        \delta X - \left[Y - Y_c(X)\right] w'(X)
        \\[0.3em]
        -\delta \ell + w(X)-\nu_{\mathrm{P}} \delta \left[Y - Y_c(X)\right]
    \end{pmatrix}.
\end{equation}

\section{Karhunen--Loève Expansion of the Random Modulus}
\label{app:kl}

The random modulus of Sec.~\ref{sec:results-exp4} has the Gaussian covariance
\begin{equation}
    \mathcal{C}(\x,\x') = \exp\left( -\frac{\|\x-\x'\|^2}{2\ell^2} \right),
\end{equation}
with $\ell=1.5$ in reference coordinates, corresponding to a correlation length of 7.5~cm. Its Karhunen--Loève modes are the eigenfunctions of the covariance operator, $\int_\Omega\mathcal{C}(\x,\x')\phi_r(\x')\,\d\x'=\lambda_r\phi_r(\x)$, ordered by decreasing eigenvalue $\lambda_r$. Since the Gaussian kernel factorizes over coordinates, each three-dimensional mode is a product of one-dimensional eigenfunctions, computed by a Nyström discretization with 64 Gauss--Legendre points. The amplitudes in Eq.~\eqref{eq:E_KL} are $a_r=\sigma\sqrt{3\lambda_r}$, with $\sigma=0.05$, where the factor 3 compensates for the variance $1/3$ of a uniform variable on $[-1,1]$. Each term of the expansion is thus a product of one-dimensional functions, and the modulus has a compact QTT representation.

Table~\ref{tab:exp4-variance} reports the fraction of the variance of the modulus captured by $M$ modes, together with that of the response $\delta$ in Eq.~\eqref{eq:exp4-response}. With $M=10$, the expansion captures 99.3\% of the material variance and keeps $E\geq0.75E_0$, a conservative bound obtained from the maxima of the modes. The response variance, in contrast, is already 99.5\% of its $M=10$ value at $M=1$: for this correlation length, the average compression is governed by the first mode.

\begin{table}[h]
    \centering
    \caption{Captured variance of the material modulus and of the response $\delta$ versus the number of modes $M$, for $N=N_\xi=8$. The response variance is normalized by its value at $M=10$.}
    \label{tab:exp4-variance}
    \begin{tabular}{l*{10}{c}}
          \toprule
          $M$ & 1 & 2 & 3 & 4 & 5 & 6 & 7 & 8 & 9 & 10 \\
          \midrule 
          Material [\%] & 67.519 & 76.407 & 85.294 & 94.182 & 95.352 & 96.522 & 97.692 & 98.227 & 98.763 & 99.298 \\
          Response [\%] & 99.515 & 99.672 & 99.805 & 99.938 & 99.955 & 99.972 & 99.980 & 99.989 & 99.994 & 100.000 \\
          \bottomrule
      \end{tabular}
\end{table}

\FloatBarrier
\bibliographystyle{elsarticle-num-names}
\bibliography{references}

\end{document}